\documentclass[a4paper,10pt,papersize]{article}

\usepackage{tikz}
\usepackage{pgfplots}

\usepackage{amsfonts,amsthm,amssymb,amsbsy,amsmath}
\usepackage{mathtools}
\usepackage{ dsfont }
\usepackage{booktabs}
\usepackage{hyperref}
\usepackage{cite,caption}
\usepackage{graphicx,hyperref}
\usepackage{algorithm}
\usepackage{algorithmicx}
\usepackage{algpseudocode}
\usepackage[varg]{txfonts}
\newcommand{\N}{\mathbb{N}}
\newcommand{\R}{\mathbb{R}}
\newcommand{\Z}{\mathbb{Z}}

\newcommand{\fz}{\frac}
\newcommand{\prz}[2]{ \frac{\partial{#1}}{\partial{#2}} }
\newcommand{\pz}{\partial}

\renewcommand{\epsilon}{\varepsilon}

\newcommand{\dps}{\displaystyle}
\renewcommand{\Omega}{\varOmega}
\renewcommand{\Gamma}{\varGamma}
\renewcommand{\Psi}{\varPsi}
\renewcommand{\Pi}{\varPi}
\newcommand{\defeq}{\coloneqq}

\newcommand{\mcL}{\mathcal{L}}
\newcommand{\sym}{{\rm sym}}
\newcommand{\uctd}[1]{ {\stackrel{\raisebox{-.1ex}{ \hskip -.1cm$\scriptscriptstyle{\bigtriangledown}$}}{ {#1} }} }

\newtheorem{Thm}{{Theorem}}

\newtheorem{Lem}{Lemma}
\newtheorem{Hyp}{Hypothesis}
\newtheorem{Rmk}{Remark}
\newtheorem{Ex}{Example}
\newtheorem{Cor}{Corollary}

\providecommand{\keywords}[1]{\textit{Keywords:} #1}

\allowdisplaybreaks[4]

\title{
Second-Order Finite Difference Approximations\\
of the Upper-Convected Time Derivative
}
\author{
D{\'e}bora O. Medeiros\thanks{Departamento de Matem{\'a}tica Aplicada e Estat{\'i}stica, Instituto de Ci{\^e}ncias Matem{\'a}ticas e de Computa\c{c}{\~a}o - ICMC Universidade de S{\~a}o Paulo - Campus de S{\~a}o Carlos
  (\texttt{deboramedeiros@usp.br}).},\ \ 
Hirofumi~Notsu\thanks{Faculty of Mathematics and Physics, Kanazawa University, Kanazawa 920-1192, Japan
  (\texttt{notsu@se.kanazawa-u.ac.jp}).},\ \ 
and\ \ 
Cassio M. Oishi\thanks{Departamento de Matem{\'a}tica e Computa\c{c}{\~a}o, Faculdade de Ci{\^e}ncias e Tecnologia Universidade Estadual Paulista J{\'u}lio de Mesquita Filho, Presidente Prudente 19060-900, SP, Brazil (\texttt{cassio.oishi@unesp.br}).}
}
\date{}
\begin{document}
\maketitle
\begin{abstract}
In this work, new finite difference schemes are presented for dealing with the upper-convected time derivative in the context of the generalized Lie derivative. The upper-convected time derivative, which is usually encountered in the constitutive equation of the popular viscoelastic models, is reformulated in order to obtain approximations of second-order in time for solving a simplified constitutive equation in one and two dimensions. The theoretical analysis of the truncation errors of the methods takes into account the linear and quadratic interpolation operators based on a Lagrangian framework. Numerical experiments illustrating the theoretical
results for the model equation defined in one and two dimensions are included. Finally, the finite difference approximations of second-order in time are also applied for solving a two-dimensional Oldroyd-B constitutive equation subjected to a prescribed velocity field at different Weissenberg numbers.
\par
\vspace{0.5em}
\noindent
\end{abstract}
%
%
%
\section{Introduction}
\label{sec:Introduction}
%
%
%
%
The solution of constitutive equations for viscoelastic fluids involves some important considerations, as for instance, the theoretical issues concerning the existence results~\cite{Chupin2018,Ervin2003,Lukakova2017,Renardy1991}, and the development of numerical schemes for solving complex fluid flows~\cite{Dellar2014,Han2000,Lee2004}. 
\par
Some forms of viscoelastic constitutive equations can be constructed considering the
upper-convected time derivative or Oldroyd derivative~\cite{Oldroyd1950}, which is defined as
\begin{equation}
	\uctd{\zeta} \defeq \frac{\partial {\zeta}}{\partial t} + \left( {u} \cdot \nabla  \right) {\zeta} - (\nabla {u}) {\zeta} - {\zeta}(\nabla {u})^\top,
	\label{der}
	\end{equation}
where $u(x,t) \in \R^d$ is the velocity field of the flow and ${\zeta}(x,t) \in \R^{d\times d}_\sym$ is a tensor to represent the non-Newtonian contribution for $d=(1,) 2, 3$.
Roughly speaking, the derivative form of~\eqref{der} is generally used for describing responses of viscoelastic fluids, as for instance, the deformation induced by the rate of strain. Therefore, the upper-convected time derivative~\eqref{der} is employed to formulate the constitutive equations of the 
most popular models, as for instance the Oldroyd-B, Phan-Tien--Tanner~(PTT), Giesekus, etc~\cite{MalPruSkrSul-2018,RenardyBook}. 
\par
In particular, we are interested in the numerical approximations for model equations based on the classical differential constitutive equation for the Oldroyd-B fluid in a dimensionless form:
	\begin{equation}
	{\zeta} + Wi \, \uctd{\zeta} = 2\left( 1 - \beta \right) {D}({u}),
	\label{eq:OB}
	\end{equation}
where ${D}({u}) = {[\nabla {u} + (\nabla {u})^\top]/{2}}$ is the strain-rate tensor, and the non-dimensional positive parameters $Wi$ and $\beta$ are respectively the Weissenberg number and the viscosity ratio~$(\beta \in (0,1))$.
\par
The Weissenberg number \cite{White:1964} is a parameter related to the memory of the fluid, i.e., for a viscoelastic material, the $Wi$ is a dimensionless number which can represent the relaxation time of the fluid. From a rheological point-of-view, the Weissenberg number can be interpreted as a number which can be used to measure the competition between elastic and viscous forces present in the concept of the viscoelasticity. A naive form to interprete the mathematical effect of this non-dimensional number is considering if $Wi=0$ in Eq. (1.2), and in this case, the stress, represented here by $\zeta$, is given by an explicit relation with the strain-rate tensor $D(u)$. Otherwise, for $Wi\neq 0$, the relation between the stress and the velocity gradient (rate-of-strain) can be modeled by a differential model, as for instance Eq. (1.2). Notice that increasing the value of the Weissenberg number in Eq. (1.2), the convected time derivative assumes a more significant effect in the equation, and therefore, the numerical treatment of this term needs to be improved in order to obtain a correct approximation of the solution. More details concerning the effect of the Weissenberg number on the partial differential equations whose describe viscoelastic fluid flows can be found in the works of Renardy \cite{RenardyAnnual,RenardyBook}.
\par
From a numerical point of view, in order to preserve the stability of the solutions, Eulerian frameworks for solving equation~\eqref{eq:OB} need to apply a high-order spatial discretization for treating the convective terms in~\eqref{der}. Generally, the methods for dealing with convection-dominant terms of the upper-convected time derivative are based on the explicit and implicit upwind methodologies~\cite{Baba1981,Harten1984,Schlichting2017}. Considering explicit upwind strategies, many numerical approaches have been proposed in the literature for solving constitutive equations of viscoelastic models based on Eq.~\eqref{eq:OB}, e.g. the Eulerian schemes using Finite-Element~(FE) \cite{Castillo,Fortin,Hulsen,Sandri}, Finite-Volume~(FV) \cite{Alves2003,Darwish,Oliveira,Pimenta}, Finite-Difference~(FD)~\cite{Franca,Martins2015,Tome}, etc. It is worth to notice that the main drawback of the explicit upwind schemes is the severe time step limitations, and the application of implicit time integrators has been used for developing more robust frameworks \cite{Breuss2006,Schlichting2017,Yee1996}, where a typical example is the so-called CFL condition. However, the construction of fully implicit upwind algorithms is complex resulting in general in high-cost computational schemes due to the solution of large systems. An additional drawback of implicit upwind schemes for solving convection-dominant problems is the excessive numerical diffusion. 
\par
In a different framework, Lagrangian methods combined with the method of characteristics~\cite{Benitez2012,BermejoGalanSaavedra2012,Douglas1982,NT-2016-M2AN,Sul-1988} for solving viscoelastic fluid flows have been proposed by~\cite{Baranger1997,Basombrio1991,Hadj1990,Notsu2015,LMNT-Peterlin_Oseen_Part_I,LMNT-Peterlin_Oseen_Part_II,Machmoum2001}. In these schemes, the Eulerian discretization of the convective term in~\eqref{der}, i.e., $(u\cdot\nabla) \zeta$, is avoided by using a Lagrangian discretization of the material derivative, i.e., $\pz\zeta/\pz t + (u\cdot\nabla) \zeta$, with the idea of the method of characteristics.
The idea is to consider the trajectory of a fluid particle and discretize the material derivative along the trajectory.
Since it is natural from a physical viewpoint and such Lagrangian schemes have advantages, e.g., the symmetry of resulting coefficient matrices of the system of linear equations in the implicit framework, no artificial parameters and no need of the so-called CFL condition, they are useful for flow problems appearing in the field of scientific computing.
\par
A different approach for avoiding numerical instabilities and to obtain accurate solutions of Eq.~\eqref{eq:OB} is mathematically rooted on the concept of the generalized Lie derivatives~(GLD)~\cite{Lee2006,LeeXuZhang2011,LeeBook2011} which modifies the definition of Eq.~\eqref{der}. In particular, this elegant methodology was firstly presented by Lee and Xu~\cite{Lee2006} (see also a similar idea proposed in~\cite{Petera2002}). In that pioneer work, the authors reformulated Eq.~\eqref{eq:OB} using some mathematical properties to define generalized Ricatti equations in terms of GLD. In summary, the upper-convected time derivative~\eqref{der} was re-written using the concept of the transition matrix. This idea was adopted in the context of the finite element discretization in Lee et al. \cite{LeeBook2011} to numerically solve the Poiseuille flow between two parallel plates around a cylinder while in~\cite{Lee2006} the authors presented theoretical results concerning the discretized version of the formulation proposed in~\cite{LeeBook2011}.
\par
In spite of the good stability properties observed in the numerical results and the sophisticated theoretical analysis of the works in~\cite{Lee2006,LeeBook2011}, to the best knowledge of the authors, the application of the GLD for solving equations in the form of~\eqref{eq:OB} is limited for finite element discretization resulting in schemes of (mainly) first-order in time.
In~\cite{Lee2006}, two finite element schemes of second-order in time are presented based on the Crank--Nicolson or the Adams--Bashforth method along the trajectory of fluid particle. There are, however, no truncation error analysis of second-order in time and no numerical results yet, while numerical results by a GLD-based finite element scheme of first-order in time are given in~\cite{LeeBook2011}.
Therefore, main contributions of this work can be summarized as follows: i)~the combination of the GLD strategy with the method of characteristics to develop temporal second-order finite difference schemes for treating the upper-convected time derivative~\eqref{der}, and ii)~the application of simple stable algorithms avoiding the need to solve large systems as commonly occur for implicit upwind schemes.
\par
In this paper, we present finite difference approximations of the upper-convected time derivative~\eqref{der} based on GLD, and apply them to simple models.
The approximations are of second-order in time, where the truncation error of second-order in time is proved in Theorem~\ref{thm:sec_order}, and a practical form is given in Corollary~\ref{cor:sec_order}.
To the best knowledge of the authors, it is noted that the form, cf.~\eqref{eq:sec_order_n}, in the corollary is new and that there are no proofs of truncation error of second-order in time for time-discretized approximations using GLD-approach.
Combining the approximation with the (bi)linear~($p=1$) and (bi)quadratic~($p=2$) Lagrange interpolations, we present full discretizations of the upper-convective time derivative of second-order in time and $p$-th order in space, i.e., $O(\Delta t^2+h^p)$, which are proved in Theorem~\ref{prop:sec_order_p_order}.
We present two numerical schemes for simple models in $d$-dimensional spaces~($d=1,2$), cf.~\eqref{scheme}, which are both explicit.
The difference of the schemes is the accuracy in space, i.e., one is of first-order~($p=1$) and the other is of second-order~($p=2$) in space as (bi)linear and (bi)quadratic Lagrange interpolation operators have been employed, respectively.
After the presentation of the schemes, numerical experiments for simple models in $d$-dimensional spaces~($d=1,2$) are presented.
They are consistent with the theoretical accuracies shown in Theorem~\ref{prop:sec_order_p_order}.
\par
In the case of Lagrangian finite element methods (often called Lagrange--Galerkin methods), a numerical integration is often employed in real computation for an integration of a composite function, since it is not easy to compute the integration of a composite function exactly.
In fact, a rough numerical integration may cause instability, cf.~\cite{Tab-2007,TabFuj-2006}, where a robustness of a scheme of second-order in time with a choice of $\Delta t$ depending on $h$ is discussed in the papers.
On the other hand, a quadrature-free scheme is proposed by using a mass-lumping technique in~\cite{PirTab-2010}, and schemes with the exact integration of a composite function are proposed by introducing a linear interpolation of the velocity and implemented in two-dimensional numerical experiments in~\cite{TabUch-2016-CD,Tabata2018}.
In these quadrature-free schemes, there is no discrepancy between the theory and real computation.
Besides them, to the best of our knowledge, it is still a standard technique for the integration of a composite function to employ a high-order quadrature rule, cf., e.g., \cite{BermejoSaavedra2012,ColeraCarpioBermejo2021,HNY-2016,NT-2015-JSC,NT-2016-M2AN}, whose computation cost depends mainly on the number of quadrature points.
In the end, we need to choose a suitable high-order quadrature rule by considering the computation cost and the \emph{error} depending on the (expected) solution, $\Delta t$, $h$ and so on.
In the case of Lagrangian finite difference method, however, there is no need to choose a quadrature rule as no integration is used. 
This is an advantage of the Lagrangian finite difference method, cf~\cite{NRT-2013}.
The GLD-type Lagrangian finite difference schemes which will be presented in this paper also have this advantage.
\par
The paper is organized as follows.
In Section~\ref{sec:Preliminaries}, basic concepts for the for flow map and the upper-convected time derivative in the framework of the generalized Lie derivative and a simple model to be dealt in this paper are introduced.
In Section~\ref{sec:FD_discretizations}, finite difference discretizations  of the upper-convected time derivative are presented, where truncation errors are proved.
In Section~\ref{nm}, GLD-type numerical schemes of second-order in time and $p$-th order in space for the simple model and their algorithms are presented.
In Section~\ref{numerics}, numerical results by our schemes are presented to see the experimental orders of convergence.
In Section~\ref{sec:conclusions}, conclusions are given.
In Appendix, properties of GLD introduced in Section~\ref{sec:Preliminaries} are proved, and the main algorithms of the work are described in details.
%
%
%
%
\section{Preliminaries}
\label{sec:Preliminaries}
%
%
%
%
In this section, we present some basic concepts concerning the flow map and the ideas of the generalized Lie derivatives. For these purposes, we need to consider some mathematical statements.
\par
Let $\Omega \subset \R^d~(d=1, 2, 3)$ be a bounded domain and $T$ be a positive constant.
Let $u:\Omega\times (0,T) \to \R^d$ be a given velocity with the following hypothesis:
\begin{Hyp}\label{hyp:u}
	The velocity $u$ is sufficiently smooth and satisfies $u_{|\pz\Omega}=0$.
\end{Hyp}
Let $\Delta t>0$ be a time increment, $N_T \defeq \lfloor T/\Delta t \rfloor$ the total number of time steps, and $t^n \coloneqq n\Delta t~(n\in\Z)$.
For a function~$f$ defined in $\Omega\times (0,T)$, let $f^n \defeq f(\cdot,t^n)$ be the function at $n$-th time step.
We define two mappings $X_1, \tilde{X}_1: \Omega\times (0,T)\to \R^d$ by
\[
X_1(x,t)  \defeq x-\Delta t\, u(x,t), \qquad \tilde{X}_1(x,t)  \defeq x-2\Delta t\, u(x,t),
\]
which are upwind points of $x$ with respect to $u(x,t)$.
{We introduce a symbol ``$\circ$'' to represent a composition of functions defined by}
\[
(g \circ X_1^n) (x) \coloneqq g( X_1^n (x) ),
\]
for a function~$g$ defined in~$\Omega$, where $X_1^n(x) = X_1(x,t^n) = x - \Delta t\,u^n(x)$.
We prepare a hypothesis for $\Delta t$:
\begin{Hyp}\label{hyp:dt}
	The time increment $\Delta t$ satisfies
	$\dps \Delta t |u|_{C^0([0,T];W^{1,\infty}(\Omega)^d)} \le 1/8$.
\end{Hyp}
\begin{Rmk}\label{rmk:upwind_points}
	Hypotheses~\ref{hyp:u} and~\ref{hyp:dt} ensure that $X_1(\Omega, t) = \tilde{X}_1(\Omega, t) = \Omega$, and that {\rm Jacobian}s of the mappings~$X_1(\cdot,t)$ and $\tilde{X}_1(\cdot,t)$ are greater than or equal to $1/2$, for $t\in [0,T]$, cf.~\cite{RuiTab-2002,Tabata2018}.
	We note that Hypothesis~\ref{hyp:dt} has no relation with the so-called CFL condition as any spatial mesh size is not included in it.
\end{Rmk}
%
%
%
%
\subsection{Lagrangian framework and the generalized Lie derivative}
%
%
%
%
\par
\begin{figure}[!htbp]
		\centering
		\begin{tikzpicture}[thick,scale=3, every node/.style={scale=1.7}]
		\draw[thick,->] (-1.0,-1.0) -- (1.1,-1.0);
		\draw[thick,->] (-1.0,-1.0) -- (-1.0,1.1);
		
		\draw[dashed] (-1.0,-0.25) -- (1.0,-0.25);
		\draw[dashed] (-1.0,0.5) -- (1.0,0.5);
		
		\draw[blue]   (-0.75,-1.0) to[out=90,in=-100] (-0.25,1.0);
		\draw[blue]   (-0.125,-1.0) to[out=90,in=-100] (0.375,1.0);
		\draw[blue]   (0.5,-1.0) to[out=90,in=-100] (1.0,1.0);
		
		\filldraw[fill] (-0.125,-1.0) circle (0.02cm)
		(-0.01,-0.25) circle (0.02cm)
		(0.25,0.5) circle (0.02cm)
		(-0.75,-1.0) circle (0.02cm)  
		(-0.64,-0.25) circle (0.02cm) 
		(-0.375,0.5) circle (0.02cm) 
		;
		
		\draw (-0.15,-1.2) node{$(x,t)$};
		\draw (-0.75,-1.2) node{\scriptsize $(\tilde{x},t)$}; 
		\draw (1.25,-1.2) node{$\mathds{R}$};
		\draw (-1.15,-1.0) node{$t$};
		\draw (-1.15,-0.25) node{$s$};
		\draw (-1.15,0.5) node{$\tilde{t}$};
		\draw (-1.25,1.2) node{$time$};
		\draw (-0.6,-0.15) node{\scriptsize $X(\tilde{x},t;s)$}; 
		\draw (0.425,-0.15) node{$X(x,t;s)$};
		\draw (-0.35,0.6) node{\scriptsize $X(\tilde{x},t;\tilde{t})$}; 
		\draw (0.66,0.6) node{$X(x,t;\tilde{t})$};

		\end{tikzpicture}
		\caption{Sketch of the flow map for $X(x,t; s)$.}
		\label{fig:flowmap}
	\end{figure}
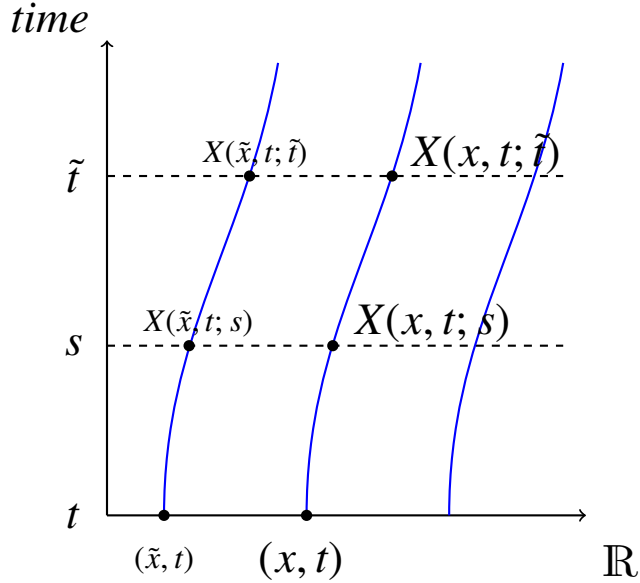

For a fixed $(x,t)\in \bar\Omega\times [0,T]$, let $X(x,t; s) \in \R^d$ be a solution of the following ordinary differential equation with an initial condition:
	\begin{subequations}\label{eqns:ode}
		\begin{align}
		\prz{}{s} X(x,t; s) & = u(X(x,t; s), s),\quad s \in (0, T),\label{eq:ode1}\\
		X(x,t; t) & = x, \label{eq:ode2}
		\end{align}
	\end{subequations}
	for $(x,t) \in \Omega\times (0,T)$.
	Physically, $X(x,t; s)$ gives the position of fluid particle at time~$s$ whose position at time~$t$ is $x$. It is known as a flow map and an illustration of this concept can be seen in Fig.~\ref{fig:flowmap}.
	\par
	For $(x,t) \in \Omega\times (0,T)$, let us introduce a matrix valued function~$L(x,t; \cdot, \cdot): (0,T)\times (0,T) \to \R^{d\times d}$ defined by
	\begin{eqnarray}\label{def:L}
	L_{ij} (x,t; t_1, t_2) \defeq \Bigl[ \prz{}{z_j} X_i(z,t_1; t_2) \Bigr]_{{\dps |}z=X(x,t; t_1)},\quad i,j = 1, \ldots, d,
	\end{eqnarray}
	which is the so-called deformation gradient.
	It is known that the function~$L$ has the following properties:
	\begin{subequations}\label{eqns:L}
		\begin{align}
		L(x,t;t_1,t_2) L(x,t;t_2,t_1) & = L(x,t;t_1,t_1) = I,
		\label{eqns:L1} \\
		\prz{}{s} L(x,t;t_1,s) & = (\nabla u) \bigl( X(x,t; s), s \bigr) L(x,t;t_1,s), 
		\label{eqns:L2} \\
		\prz{}{s} L(x,t;s,t_1) & = - L(x,t;s,t_1) (\nabla u) \bigl( X(x,t; s), s \bigr), 
		\label{eqns:L3}
		\end{align}
	\end{subequations}
	for $t_1, t_2 \in [0, T]$, where $I \in \R^{d\times d}_\sym$ is the identity matrix.
	Although the proofs can be found in, e.g., \cite{LeeBook2011}, we give the proofs again in Appendix~\ref{A.subsec:L} under the assumption of unique existence of smooth regular~$L$.
	\par
	Let $D/Dt$ be the material derivation defined by
	\[
	\fz{D}{Dt} \defeq \prz{}{t} + u\cdot\nabla.
	\]
	For a function~$\zeta:\Omega\times(0,T)\to\R^{d\times d}$, it is well-known that the material derivative of~$\zeta$ can be written as
	\begin{equation}
	\fz{D\zeta}{Dt} (x,t) = \Bigl[ \prz{\zeta}{t} + (u\cdot\nabla) \zeta \Bigr] (x,t) = \frac{\pz}{\pz s} \zeta \bigl( X(x,t; s), s \bigr)_{{\dps |}s=t}.
	\end{equation}
	Here, we define the so-called generalized Lie derivative~$\mcL_u\zeta$ by
	\begin{align}
	(\mcL_u\zeta) \bigl( X(x,t; s), s \bigr) 
	& \defeq
	L(x,t;t,s) \prz{}{s}\Bigl[ L(x,t;s,t) \zeta \bigl( X(x,t; s), s \bigr) L(x,t;s,t)^\top \Bigr] L(x,t;t,s)^\top.
	\label{def:GLD}
	\end{align}
	From~\eqref{eqns:L}, the upper-convected time derivative can be rewritten by using $\mcL_u\zeta$, i.e.,
	\begin{equation}
	\uctd{\zeta}(x,t) 
	= (\mcL_u\zeta) (x,t) = (\mcL_u\zeta)\bigl( X(x,t; s), s \bigr)_{|s=t},
	\label{eq:UCM_GLD}
	\end{equation}
	which is shown in Appendix~\ref{A.subsec:UCM_GLD}. 
%
%
%
%
\subsection{The model equation}
\label{sce}
%
%
%
%
Based on the above description, we consider a simplified model equation in order to present the application of finite difference schemes for dealing with the generalized Lie derivative. Particularly, based on the Oldroyd-B constitutive equation (\ref{eq:OB}), the problem is to find $\zeta: \Omega\times (0, T) \to \R^{d\times d}_\sym$ such that
	\begin{subequations}\label{eqns:prob}
		\begin{eqnarray}
		\uctd{\zeta} & = & F \quad\ \mbox{in}\ \Omega\times (0, T), 
		\label{eq:prob1}\\
		\zeta & = & \zeta_\mathrm{in} \quad \mbox{on}\ \Gamma_\mathrm{in}\times (0, T), 
		\label{eq:prob2}\\
		\zeta & = & \zeta^0 \quad\ \mbox{in}\ \Omega,\ \mbox{at}\ t=0, 
		\end{eqnarray}
	\end{subequations}
	where $\Gamma_\mathrm{in}$ is an inflow boundary defined by $\Gamma_\mathrm{in} \defeq \{ x\in\pz\Omega;\ u(x,t)\cdot n(x) < 0 \}$ for the outward unit normal vector~$n: \pz\Omega\to\R^d$, and $F: \Omega\times (0,T) \to \R^{d\times d}_\sym$, $\zeta_\mathrm{in}: \Gamma_{\rm in}\times (0, T) \to \R^{d\times d}_\sym$ and $\zeta^0: \Omega \to \R^{d\times d}_\sym$ are given functions.
	\begin{Rmk}
		(i)~From~\eqref{eq:UCM_GLD}, Eq.~\eqref{eq:prob1} can be reformulated using the generalized {\rm Lie} derivative resulting in:
	\begin{equation}
	\mcL_u\zeta
	= F \quad \mbox{in}\ \Omega\times (0, T).
	\label{eq:prob1nf}
	\end{equation}
	(ii)~In general, the inflow boundary~$\Gamma_\mathrm{in}$ depends on time~$t$, i.e., $\Gamma_\mathrm{in} = \Gamma_\mathrm{in}(t)$, while $\Gamma_\mathrm{in}$ is the empty set under Hypothesis~\ref{hyp:u}.
	Throughout the paper, we deal with the inflow boundary~$\Gamma_\mathrm{in}$ independent of time~$t~(\in (0,T))$.
	\end{Rmk}
%
%
%
%
%
%
%
%
\section{Finite difference discretizations}
\label{sec:FD_discretizations}
%
%
%
%
In this section, we present descriptions concerning the spatial and temporal discretizations. The main results related to the numerical analysis of the schemes are also described in details.
%
%
%
%
\subsection{Space discretizations and interpolation operators}
%
%
%
%
In this subsection, we introduce spatial discretizations and interpolation operators in one- and two-dimensions.
Before starting them, for an integer~$i$ and a positive number~$\delta$, we prepare two functions $\eta_i^{(1)}(\,\cdot\,; \delta)$ and $\eta_i^{(2)}(\,\cdot\,; \delta):~\R \to \R$.
The former, $\eta_i^{(1)}(\,\cdot\,; \delta)$, is defined by
\begin{align*}
\dps
\eta_i^{(1)}(s; \delta) & \defeq
\left\{
\begin{aligned}
& \fz{s}{\delta} -i + 1 && \bigl( s\in [(i-1)\delta, i\delta ) \bigr),\\
& i+1 - \fz{s}{\delta} && \bigl( s\in [i\delta, (i+1)\delta] \bigr),\\
& \ \ \ \ \   0 && \bigl( {\rm otherwise} \bigr),
\end{aligned}
\right.
\end{align*}
and the latter, $\eta_i^{(2)}(\,\cdot\,;\delta)$, is defined by
\begin{align*}
\intertext{(i) $i$ : even number}
\dps
\eta_i^{(2)}(s;\delta) & \defeq
\left\{
\begin{aligned}
& \Bigl( \fz{s}{\delta} - i + 1 \Bigr) \Bigl( \fz{s}{2\delta} - \fz{i}{2} + 1 \Bigr) && \bigl( s\in [(i-2)\delta, i\delta ) \bigr),\\
&\Bigl( i+1 - \fz{s}{\delta} \Bigr) \Bigl( \fz{i}{2} + 1 - \fz{s}{2\delta} \Bigr)  && \bigl( s\in [i\delta, (i+2)\delta] \bigr),\\
& \qquad\qquad\ \ \  0 && \bigl( {\rm otherwise} \bigr),
\end{aligned}
\right.
\intertext{(ii) $i$ : odd number}
\eta_i^{(2)}(s;\delta) & \defeq
\left\{
\begin{aligned}
& \Bigl( \fz{s}{\delta} - i + 1\Bigr) \Bigl( i+1 - \fz{s}{\delta} \Bigr) && \bigl( s\in [(i-1)\delta, (i+1)\delta] \bigr),\\
& \qquad\quad\quad \ \ \  0 && \bigl( {\rm otherwise} \bigr).
\end{aligned}
\right.
\end{align*}
The functions~$\eta_i^{(1)}(\,\cdot\,; \delta)$ and~$\eta_i^{(2)}(\,\cdot\,; \delta)$ are used below for the definitions of (bi)linear and (bi)quadratic interpolation operators $\Pi_h^{(1)}$ and $\Pi_h^{(2)}$, respectively.
%
%
%
%
\subsubsection{One-dimensional case $(d=1)$}
%
%
%
%
Initially, we consider one spatial dimension, i.e., $d = 1$.
For the sake of simplicity, we assume $\Omega = (0, a)$ for a positive number~$a$. 
Let $N\in\N$ be a number, $h\defeq a/N$ a mesh size, and $x_i \defeq ih~(i\in\Z)$ lattice points.
We define a set of lattice points~$\bar\Omega_h$ and a discrete function space~$V_h$ restrict to the number $N$, by
\begin{align*}
\bar\Omega_h & 
\defeq \{ x_i \in \bar\Omega;\ i = 0,\ldots,N \} \ (\subset \bar\Omega \subset \R^d = \R), \\
V_h & \defeq \{ v_h: \bar\Omega_h \to \R^{d\times d}_\sym \} = \{ v_h: \bar\Omega_h \to \R \}.
\end{align*}
We introduce a set of basis functions~$\{\varphi_i^{(1)}: \bar\Omega \to \R; \ i=0, \ldots, N \}$ defined by
\[
\varphi_i^{(1)}(x) \defeq \eta_i^{(1)}(x; h),
\qquad
i=0,\ldots,N.
\]
The functions~$\varphi_0^{(1)}$ and $\varphi_N^{(1)}$ are simplified to
\begin{align*}
\dps
\varphi_0^{(1)}(x) & \defeq
\left\{
\begin{aligned}
& 1 - \fz{x}{h} && \bigl( x\in [x_0, x_1] \bigr) \\
& \ \ \   0 && \bigl( {\rm otherwise} \bigr)
\end{aligned}
\right.
=
\left\{
\begin{aligned}
&\fz{x_1-x}{h} && \bigl( x\in [ x_0, x_1] \bigr) \\
& \ \ \ \   0 && \bigl( {\rm otherwise} \bigr)
\end{aligned}
\right.,
\\
\varphi_N^{(1)}(x) & \defeq
\left\{
\begin{aligned}
& \fz{x}{h} - N + 1 && \bigl( x\in [x_{N-1}, x_N ] \bigr) \\
& \ \ \ \ \ \   0 && \bigl( {\rm otherwise} \bigr)
\end{aligned}
\right.
=
\left\{
\begin{aligned}
&\fz{x-x_{N-1}}{h} && \bigl( x\in [ x_{N-1}, x_N ] \bigr) \\
& \ \ \ \ \ \   0 && \bigl( {\rm otherwise} \bigr)
\end{aligned}
\right.,
\end{align*}
as defined in $\bar\Omega = [x_0, x_N] = [0, a]$.
Let $\Pi_h^{(1)}: V_h \to C^0(\bar\Omega)$ be the linear interpolation operator defined by
\[
\bigl( \Pi_h^{(1)} v_h \bigr) (x) \defeq \sum_{i=0}^N v_h(x_i)\varphi_i^{(1)}(x).
\]
\par
We describe the ideas for using a quadratic interpolation. Let $N\in\N$ be an even number, and $M\defeq N/2\in\N$. For the definition of the quadratic interpolation operator $\Pi_h^{(2)}$, we define a set of basis functions~$\{\varphi_i^{(2)}: \bar\Omega \to \R; \ i=0, \ldots, N \}$ by
\[
\varphi_i^{(2)}(x) \defeq \eta_i^{(2)}(x; h),
\qquad
i=0,\ldots,N,
\]
where $\varphi_0^{(2)}$ and $\varphi_N^{(2)} \ (=\varphi_{2M}^{(2)})$ are reduced to
\begin{align*}
\dps
\varphi_0^{(2)}(x) & =
\left\{
\begin{aligned}
&\fz{(x_1-x)(x_2-x)}{2h^2} && \bigl( x\in [x_0, x_2] \bigr),\\
& \qquad\quad \  0 && \bigl( {\rm otherwise} \bigr),
\end{aligned}
\right.
\\
\varphi_N^{(2)}(x) & =
\left\{
\begin{aligned}
&\fz{(x-x_{N-1})(x-x_{N-2})}{2h^2} && \bigl( x\in [x_{N-2}, x_N ] \bigr),\\
& \qquad\quad\ \ \ \ \  0 && \bigl( {\rm otherwise} \bigr).
\end{aligned}
\right.
\end{align*}
Let $\Pi_h^{(2)}: V_h \to C^0(\bar\Omega)$ be the quadratic interpolation operator defined by
\[
\bigl( \Pi_h^{(2)} v_h \bigr) (x) \defeq \sum_{i=0}^N v_h(x_i)\varphi_i^{(2)}(x).
\]
\begin{Rmk}\label{rmk:upwind_cell_1d}
	For $\alpha$, $\beta \in \R~(\alpha<\beta)$, and $N_0\in\N$ with $\delta_0 = \delta_0(\alpha,\beta,N_0) \defeq (\beta-\alpha)/N_0 > 0$, let $\mathcal{I}(\cdot; \alpha, \beta, N_0) : \R \to \{0,\ldots,N_0\}$ be an integer-valued index indicator function defined by
	\begin{align}
	\label{def:index_func}
	\mathcal{I} (s; \alpha, \beta, N_0) \defeq
	\left\{
	\begin{aligned}
	& \left\lfloor \fz{s-\alpha}{\delta_0} \right\rfloor && \bigl( s \in (\alpha,\beta) \bigr), \\
	& \quad\ \ 0 && (s \le \alpha), \\
	& \quad\ N_0 && (s \ge \beta).
	\end{aligned}
	\right.
	\end{align}
	We note that the integer~$i_0 = \mathcal{I}(s; \alpha, \beta, N_0)$ satisfies $i_0\delta_0+\alpha \le s < (i_0+1)\delta_0+\alpha$ for $s\in (\alpha, \beta)$, and that, for an even number~$N_0$ with $M_0 = N_0/2 \in\N$, the integer~$k_0 = \mathcal{I}(s; \alpha, \beta, M_0)$ satisfies $2k_0\delta_0+\alpha \le s < 2(k_0+1)\delta_0+\alpha$ for $s\in (\alpha, \beta)$ as $\delta_0 = (\beta - \alpha)/M_0 = 2 (\beta - \alpha)/N_0$. 
	\par
	For $d=1$, we introduce two notations of intervals,
	\begin{align*}
	&&&& K_{i+1/2}^{(1)} & \defeq [x_i, x_{i+1}], & i & \in\{0,\ldots,N-1\}, &&&& \\
	&&&& K_{2k+1}^{(2)} & \defeq [x_{2k}, x_{2k+2}], & k & \in\{0,\ldots,M-1\}, &&&&
	\end{align*}
	whose measures are $h$ and $2h$, respectively.
	Let $x \in \R$ be given arbitrarily.
	Then, the following are practically useful in computation: \smallskip\\
	$(i)$~Let $i_0 \defeq \mathcal{I}(x; 0, a, N)\in\{0,\ldots,N\}$.
	When $x\in\Omega$, the integer $i_0$ satisfies $x\in K_{i_0+1/2}^{(1)} = [x_{i_0}, x_{{i_0}+1}]$, and we have two-points representation of~$( \Pi_h^{(1)} v_h) (x)$,
	\begin{equation}
	\bigl( \Pi_h^{(1)} v_h \bigr) (x) = v_{i_0} \varphi_{i_0}^{(1)}(x) + v_{i_0+1} \varphi_{{i_0}+1}^{(1)}(x),
	\label{int1}
	\end{equation}
	where we have used a notation~$v_i = v_h(x_i)$. \smallskip\\
	$(ii)$~Let $k_0 \defeq \mathcal{I}(x; 0, a, M)\in\{0,\ldots,M\}$. 
	When $x\in\Omega$, the integer $k_0$ satisfies $x\in K_{2k_0+1}^{(2)} = [x_{2k_0}, x_{2{k_0}+2}]$, and we have three-points representation of~$( \Pi_h^{(2)} v_h) (x)$,
	\begin{equation}
	\bigl( \Pi_h^{(2)} v_h \bigr) (x) = v_{2k_0} \varphi_{2k_0}^{(2)}(x) + v_{2k_0+1} \varphi_{2k_0+1}^{(2)}(x) + v_{2k_0+2} \varphi_{2k_0+2}^{(2)}(x),
	\label{int2}
	\end{equation}
	for $v_i = v_h(x_i)$. \smallskip\\
	$(iii)$~If the value~$(\Pi_h^{(p)} v_h)(x)~(p=1,2)$ is needed for $x\notin\Omega$, we can employ, instead of it, the closest end value of $v_h$, i.e., $v_0 = v_h(0)~(x \le 0)$ or $v_N = v_h(a)~(x \ge a)$, while the value~$v_0$ or~$v_N$ should be given by using~$\zeta_{\rm in}$ as $x$ corresponds to an upwind point and $x\notin\bar\Omega$ means the high possibility of existence of ``inflow'' boundary near $x$.
	The function~$\mathcal{I}(\cdot;\alpha,\beta,N)$ is, therefore, also useful for $x\notin\Omega$ in the sense that $\mathcal{I}(x;\alpha,\beta,N)$ provides the closest index of lattice point.
\end{Rmk}
%
%
%
%
%
\subsubsection{Two-dimensional case $(d=2)$}

We consider two spatial dimensions, i.e., $d=2$.
For the sake of simplicity, we assume $\Omega = (0, a_1)\times (0, a_2)$ for positive numbers~$a_1$ and $a_2$.
Let $N_i \in\N~(i=1,2)$ be numbers, $h_i\defeq a_i/N_i~(i=1,2)$ mesh sizes in $x_i$-direction, $h_{\min} \defeq \min\{h_i;\ i=1,\ldots,d\}$ and $h = h_{\max} \defeq \max\{h_i;\ i=1,\ldots,d\}$ minimum and maximum mesh sizes, and $x_{i,j} \defeq (ih_1,jh_2)^\top~(i,j\in\Z)$ lattice points.
We assume a family of meshes satisfying the next hypothesis:
\begin{Hyp}\label{hyp:mesh}
There exist positive constants $h_0$, $\gamma_1$ and $\gamma_2$ such that
\[
h \in (0,h_0], \quad \mbox{and} \quad \gamma_1 \le \fz{h}{h_{\min}} \le \gamma_2.
\]
\end{Hyp}
\begin{Rmk}
The hypothesis is set for $d=2$ essentially, as it always holds for $d=1$ with $\gamma_1 = \gamma_2 = 1$.
\end{Rmk}
\par
We define a set of lattice points~$\bar\Omega_h$ and a discrete function space~$V_h$ restrict to the numbers $N_i \in\N~(i=1,2)$, by
\begin{align*}
\bar\Omega_h & 
\defeq \{ x_{i,j} \in \bar\Omega;\ i = 0,\ldots,N_1,\ j = 0,\ldots,N_2 \}, \\
V_h &\defeq \{ v_h: \bar\Omega_h \to \R^{d\times d}_\sym \},
\end{align*}
where it is noted that $\bar\Omega_h \subset \bar\Omega \subset \R^d~(= \R^2)$.
Using~$\eta_i^{(1)}(\,\cdot\,; \delta)$, we introduce a set of basis functions~$\{\varphi_{i,j}^{(1)}: \bar\Omega \to \R; \ x_{i,j}\in\bar\Omega_h,\ i,j\in\Z \}$ defined by
\begin{displaymath}
\dps
\varphi_{i,j}^{(1)}(x) = \varphi_{i,j}^{(1)}(x_1,x_2) \defeq 
\eta_i^{(1)}(x_1;h_1) \eta_j^{(1)}(x_2;h_2).
\end{displaymath}
Let $\Pi_h^{(1)}: V_h \to C^0(\bar\Omega)$ be the bilinear interpolation operator defined by
\[
\bigl( \Pi_h^{(1)} v_h \bigr) (x) \defeq \sum_{x_{i,j}\in \bar\Omega_h} v_h(x_{i,j})\varphi_{i,j}^{(1)}(x).
\]
\par
The extension of the above interpolation using the biquadratic interpolation strategy can be defined as follows. Let $N_1, N_2 \in\N$ be even numbers, and $M_i\defeq N_i/2\in\N$ for $i=1,2$.
For the definition of the biquadratic interpolation operator $\Pi_h^{(2)}$, we introduce basis functions~$\{\varphi_{i,j}^{(2)}: \bar\Omega \to \R; \ x_{i,j}\in\bar\Omega_h \}$ defined by
\[
\varphi_{i,j}^{(2)}(x) = \varphi_{i,j}^{(2)}(x_1,x_2) \defeq \eta_i^{(2)}(x_1;h_1) \eta_j^{(2)}(x_2;h_2).
\]
Let $\Pi_h^{(2)}: V_h \to C^0(\bar\Omega)$ be the biquadratic interpolation operator defined by
\[
\bigl( \Pi_h^{(2)} v_h \bigr) (x) \defeq \sum_{x_{i,j}\in\bar\Omega_h} v_h(x_{i,j})\varphi_{i,j}^{(2)}(x).
\]
\begin{Rmk}\label{rmk:upwind_cell_2d}
For $d=2$, we introduce two notations of boxes (cells),
\begin{align*}
K_{i+1/2,j+1/2}^{(1)} & \defeq [ih_1, (i+1)h_1] \times [jh_2, (j+1) h_2], 
&
(i,j) & \in \{0,\ldots,N_1-1\} \times \{0,\ldots,N_2-1\}, \\
K_{2k+1,2l+1}^{(2)} & \defeq [2kh_1, (2k+2)h_1] \times [2l h_2, (2l+2)h_2], 
& (k,l) & \in \{0,\ldots,M_1-1\} \times \{0,\ldots,M_2-1\},
\end{align*}
whose measures are $h_1h_2$ and $4h_1h_2$, respectively.
Let $x \in \R^2$ be given arbitrarily.
Then, the following are practically useful in computation: \smallskip\\
$(i)$~Let $i_0 \defeq \mathcal{I}(x_1; 0, a_1, N_1)\in\{0,\ldots,N_1\}$ and $j_0 \defeq \mathcal{I}(x_2; 0, a_2, N_2)\in\{0,\ldots,N_2\}$.
When $x\in\Omega$, the set of integers $(i_0,j_0)$ satisfies $x\in K_{i_0+1/2,j_0+1/2}^{(1)} = [i_0h_1, ({i_0}+1)h_1] \times [j_0h_2, ({j_0}+1)h_2]$, and we have four-points representation of~$( \Pi_h^{(1)} v_h) (x)$,
\begin{align}
\label{int3}
\bigl( \Pi_h^{(1)} v_h \bigr) (x) = 
\sum_{m,n=0,1} v_{i_0+m,j_0+n} \, \varphi_{i_0+m,j_0+n}^{(1)}(x),
\end{align}
where we have used a simplified notation~$v_{i,j} = v_h(x_{i,j})$. \smallskip\\
$(ii)$~Let $k_0 \defeq \mathcal{I}(x_1; 0, a_1, M_1)\in\{0,\ldots,M_1\}$ and $l_0 \defeq \mathcal{I}(x_2; 0, a_2, M_2)\in\{0,\ldots,M_2\}$.
When $x\in\Omega$, the integer $k_0$ satisfies $x\in K_{2k_0+1}^{(2)} = [x_{2k_0}, x_{2{k_0}+2}]$, and we have nine-points representation of~$( \Pi_h^{(2)} v_h) (x)$,
\begin{equation}
\label{int4}
	\bigl( \Pi_h^{(2)} v_h \bigr) (x) = \sum_{m,n=0,1,2} v_{2k_0+m,2l_0+n} \varphi_{2k_0+m,2l_0+n}^{(2)}(x).
\end{equation}
	$(iii)$~If the value~$(\Pi_h^{(p)} v_h)(x)~(p=1,2)$ is needed for $x\notin\Omega$, we can employ, instead of it, the closest end value of $v_h$, i.e., one of the values of $v_h(x_{i,j})~(x_{i,j} \in \bar\Omega_h\cap\pz\Omega)$, while the value should be given by using~$\zeta_{\rm in}$ as $x$ corresponds to an upwind point.
\end{Rmk}
\begin{Rmk}\label{rmk:3D}
We omit the extension of the interpolation operators~$\Pi_h^{(p)}~(p=1, 2)$ to the three-dimensional case, i.e., $d=3$, since it is naturally defined by introducing basis functions~$\varphi_{i,j,k}^{(p)}(x) = \varphi_{i,j,k}^{(p)}(x_1,x_2,x_3) \defeq 
\eta_i^{(p)}(x_1;h_1) \eta_j^{(p)}(x_2;h_2) \eta_k^{(p)}(x_3;h_3)$ for $p = 1, 2$ in a similar manner.
\end{Rmk}
%
%
%
%
\subsection{Time discretization: truncation error analysis}\label{subsec:time_discretization}
%
%
%
%
%
For the velocity~$u$, let $L_1,\ \tilde{L}_1: \Omega \times (0,T) \to \R^{d\times d}$ be matrices defined by
\begin{align}
L_1 (x,t) \defeq I + \Delta t (\nabla u)(x,t),
\quad
\tilde{L}_1 (x,t) \defeq I + 2\Delta t (\nabla u)(x,t),
\label{defs:L1_tilde_L1}
\end{align}
which are approximations of $L(x,t;\, t-\Delta t,t)$ and $L(x,t;\, t-2\Delta t,t)$, respectively, cf.~Lemma~\ref{lem:L_with_k} below.
Now, we present a theorem which provides an approximation of the upper-convected time derivative of second-order in time.
\begin{Thm}\label{thm:sec_order}
	Suppose that Hypotheses~\ref{hyp:u} and~\ref{hyp:dt} hold true.
	Let $\zeta: \bar\Omega\times [0,T] \to \R^{d\times d}$ be a sufficiently smooth function.
	Then, for any $x\in\bar\Omega$ and $t\in [2\Delta t, T]$, we have
	\begin{align}
	\uctd{\zeta} (x,t) 
	= \fz{1}{2\Delta t}\Bigl[
	3\zeta(x,t) \notag
	-4 L_1 (x,t) \zeta(X_1(x,t),t-\Delta t) L_1 (x,t)^\top \\
	+ \tilde{L}_1 (x,t) \zeta(\tilde{X}_1(x,t),t-2\Delta t) \tilde{L}_1 (x,t)^\top
	\Bigr] + O(\Delta t^2).
	\label{eq:sec_order}
	\end{align}
\end{Thm}
We give the proof of Theorem~\ref{thm:sec_order} after giving a remark and preparing two lemmas.
\begin{Rmk}
(i)~Let us consider $(x,t) \in \bar\Omega\times [2\Delta t,T]$ as a fixed point and employ simple notations~$X = X(x,t;\,\cdot\,)$ and~$L(\,\cdot\,,\,\cdot\,) = L(x,t;\,\cdot\,,\,\cdot\,)$.
Then, an approximation of~$\uctd{\zeta} (x,t)$ of first-order in time is obtained as follows:
\begin{align*}
\uctd{\zeta}(x,t) & = (\mcL_u\zeta)\bigl( X(s), s \bigr)_{{\dps |}s=t} 
\qquad\qquad\qquad\qquad\qquad\qquad\qquad\qquad\qquad\quad
 \mbox{\rm (by~\eqref{eq:UCM_GLD})} \\
& = L(t,s) \prz{}{s}\Bigl[ L(s,t) \zeta \bigl( X(s), s \bigr) L(s,t)^\top \Bigr] L(t,s)^\top_{{\dps |}s=t} 
\qquad\ \, \mbox{\rm (by definition~\eqref{def:GLD})} \\
& = L(t,s) \fz{1}{\Delta t} \Bigl[ L(s,t) \zeta \bigl( X(s), s \bigr) L(s,t)^\top \\
& \qquad - L(s-\Delta t,t) \zeta \bigl( X(s-\Delta t), s-\Delta t \bigr) L(s-\Delta t,t)^\top \Bigr] L(t,s)^\top_{{\dps |}s=t} + O(\Delta t) \\
& \qquad\qquad\qquad\qquad\qquad\qquad\qquad\quad
 \mbox{\rm (by the {\rm Euler} method with respect to~$s$)} \\
& = \fz{1}{\Delta t} \Bigl[ \zeta \bigl( X(t), t \bigr) - L(t-\Delta t,t) \zeta \bigl( X(t-\Delta t), t-\Delta t \bigr) L(t-\Delta t,t)^\top \Bigr] + O(\Delta t) \\
& \qquad\qquad\qquad\qquad\qquad\qquad\qquad\qquad\quad
 \mbox{\rm (by substituting $t$ into $s$ and~\eqref{eqns:L1})} \\
& = \fz{1}{\Delta t} \Bigl[ \zeta (x,t) - L_1(x,t) \zeta \bigl( X_1(x,t), t-\Delta t \bigr) L_1(x,t)^\top \Bigr] + O(\Delta t),
\end{align*}
where the last equality holds true from the initial condition~\eqref{eq:ode2} for $X$, i.e., $X(t)=x$, and the relations,
\begin{align*}
L_1(x,t) &= L(t-\Delta t,t)+O(\Delta t^2), &
X_1(x,t) &= X(t-\Delta t) + O(\Delta t^2),
\end{align*}
which will be shown in Lemmas~\ref{lem:L_with_k} and~\ref{lem:phi_g_with_k} with $k=1$ below, respectively.
\smallskip\\
(ii)~Theorem~\ref{thm:sec_order} presents an approximation of~$\uctd{\zeta} (x,t)$ of second-order in time based on the two-step {\rm Adams--Bashforth} method, i.e., for a smooth function $f: \R\to\R$,
\[
f^\prime(t) = \fz{d}{ds}f(s)_{{\dps |}s=t} = \fz{1}{2\Delta t}[ 3f(t) - 4f(t-\Delta t) + f(t-2\Delta t)] + O(\Delta t^2),
\]
in place of the {\rm Euler} method in~(i).
\end{Rmk}
\begin{Lem}\label{lem:L_with_k}
	Suppose that Hypotheses~\ref{hyp:u} and~\ref{hyp:dt} hold true.
	Let $k = 1$ or $2$ be fixed.
	Then, for any $x\in\bar\Omega$ and $t\in [k\Delta t, T]$, we have
	\begin{equation}
	L(x,t;t-k\Delta t, t) = I + k \Delta t (\nabla u)(x,t) + \fz{(k\Delta t)^2}{2} U(x,t) + O(\Delta t^3),
	\label{eq:L_with_k}
	\end{equation}
	where $U: \Omega\times (0, T) \to \R^{d\times d}$ is a function defined by
	\[
	U \defeq (\nabla u)^2 - \fz{D (\nabla u)}{Dt}.
	\]
\end{Lem}
\begin{proof}
	From the Taylor expansion, we have
	\begin{align}
	& L(x,t; t-k\Delta t, t) 
	= L(x,t; s-k\Delta t, t)_{|s=t} \notag\\
	& = \Bigl[ L(x,t; s, t) - k\Delta t \prz{}{s} L(x,t; s, t) + \fz{(k\Delta t)^2}{2} \prz{^2}{s^2} L(x,t; s, t) \Bigr]_{{\dps |}s=t} + O(\Delta t^3) \notag\\
	& = \Bigl[ L(x,t; s, t) - k\Delta t \bigl[ - L(x,t; s, t) (\nabla u)\bigl( X(x,t; s), s \bigr) \bigr] + \notag\\
	& \quad + \fz{(k\Delta t)^2}{2} \prz{}{s} \bigl[ - L(x,t; s, t) (\nabla u)\bigl( X(x,t; s), s \bigr) \bigr] \Bigr]_{{\dps |}s=t} + O(\Delta t^3) \qquad \mbox{(by~\eqref{eqns:L3})}\notag\\
	& = I + k\Delta t (\nabla u) (x,t) 
	- \fz{(k\Delta t)^2}{2} \prz{}{s} \bigl[ L(x,t; s, t) (\nabla u)\bigl( X(x,t; s), s \bigr) \bigr]_{{\dps |}s=t} + O(\Delta t^3).
	\label{eq:proof_lem_L_with_k_1}
	\end{align}
	We evaluate $\prz{}{s} \bigl[ L(x,t; s, t) (\nabla u)\bigl( X(x,t; s), s \bigr) \bigr]_{|s=t}$ as follows:
	\begin{align}
	& \prz{}{s} \bigl[ L(x,t; s, t) (\nabla u)\bigl( X(x,t; s), s \bigr) \bigr]_{{\dps |}s=t} \notag\\
	& = \Bigl[
	\Bigl( \prz{}{s} L(x,t; s, t) \Bigr) (\nabla u)\bigl( X(x,t; s), s \bigr) 
	+ L(x,t; s, t) \Bigl( \prz{}{s} (\nabla u)\bigl( X(x,t; s), s \bigr) \Bigr)
	\Bigr]_{{\dps |}s=t} \notag\\
	& = \Bigl[
	- L(x,t; s, t) (\nabla u)^2 \bigl( X(x,t; s), s \bigr) 
	+ L(x,t; s, t) \fz{D (\nabla u)}{Dt} \bigl( X(x,t; s), s \bigr)
	\Bigr]_{{\dps |}s=t} \notag\\
	& = 
	- (\nabla u)^2 (x,t)
	+ \fz{D (\nabla u)}{Dt} (x,t)
	= - U(x,t).
	\label{eq:proof_lem_L_with_k_2}
	\end{align}
	Combining~\eqref{eq:proof_lem_L_with_k_2} with~\eqref{eq:proof_lem_L_with_k_1}, we obtain~\eqref{eq:L_with_k}.
\end{proof}
\begin{Lem}\label{lem:phi_g_with_k}
	Suppose that Hypotheses~\ref{hyp:u} and~\ref{hyp:dt} hold true.
	Let $k = 1$ or $2$ be fixed.
	Then, for any $x\in\bar\Omega$ and $t\in [k\Delta t, T]$, we have the following:\\
	(i)~It holds that
	\[
	X(x,t; t-k\Delta t) = x - k \Delta t u(x,t) + \fz{(k\Delta t)^2}{2}\fz{Du}{Dt} (x,t) + O(\Delta t^3).
	\]
	(ii)~Let $\zeta: \Omega\times (0,T) \to \R^{d\times d}$ be a sufficiently smooth function.
	It holds that
	\[
	\zeta \bigl( X(x,t; t-k\Delta t), t-k\Delta t \bigr)
	 = \zeta \bigl( x - k \Delta t u(x,t), t-k\Delta t \bigr) + \fz{(k\Delta t)^2}{2} Z(x,t) + O(\Delta t^3),
	\]
	where $Z: \Omega\times (0,T) \to \R^{d\times d}$ is a function defined by
	\[
	Z \defeq \Bigl( \fz{Du}{Dt} \cdot \nabla \Bigr) \zeta.
	\]
\end{Lem}
\begin{proof}
	We prove~(i).
Recalling that $X(x,t;s)$ is a solution to~\eqref{eqns:ode} and noting that the following identity,
\[
X(x,t; t-k\Delta t) = x - \int_{t-k\Delta t}^t u\bigl( X(x,t; s), s \bigr) ds,
\]
holds true, we have
	\begin{align*}
	& X(x,t; t-k\Delta t) - [ x - k\Delta t u(x,t)]  \\
	& = x - \int_{t-k\Delta t}^t u\bigl( X(x,t; s), s \bigr) ds - \Bigl[ x - \int_{t-k\Delta t}^t u\bigl( X(x,t; t), t \bigr) ds \Bigr] \\
	& = \int_{t-k\Delta t}^t \Bigl[ u \bigl( X(x,t; t), t \bigr) - u\bigl( X(x,t; s), s \bigr) \Bigr] ds 
	= \int_{t-k\Delta t}^t ds \Bigl[ u \bigl( X(x,t; s_1), s_1 \bigr) \Bigr]_{s_1=s}^t \\
	& = \int_{t-k\Delta t}^t ds \int_s^t \fz{Du}{Dt} \bigl( X(x,t; s_1), s_1 \bigr) ds_1 
	= \int_{t-k\Delta t}^t ds \int_s^t \Bigl( \fz{Du}{Dt} (x,t) + O(\Delta t) \Bigr) ds_1 \\
	& = \fz{(k\Delta t)^2}{2} \fz{Du}{Dt} (x,t) + O(\Delta t^3),
	\end{align*}
	which completes the proof of~(i).
	\par
	We prove~(ii).
	From (i) and the Taylor expansion, we have
	\begin{align*}
	& \zeta \bigl( X(x,t; t-k\Delta t), t-k\Delta t \bigr) \\
	& = \zeta \biggl( x - k \Delta t u(x,t) + \fz{(k\Delta t)^2}{2}\fz{Du}{Dt} (x,t), t-k\Delta t \biggr) + O(\Delta t^3) \\
	& = \zeta \bigl( x - k \Delta t u(x,t), t-k\Delta t \bigr) + \fz{(k\Delta t)^2}{2} \Bigl[ \Bigl( \fz{Du}{Dt} (x,t) \cdot \nabla \Bigr) \zeta \Bigr] \bigl( x - k \Delta t u(x,t), t-k\Delta t \bigr) + O(\Delta t^3) \\
	& = \zeta \bigl( x - k \Delta t u(x,t), t-k\Delta t \bigr) + \fz{(k\Delta t)^2}{2} Z (x,t) + O(\Delta t^3),
	\end{align*}
	where we have used the relation,
	\[
	\Bigl[ \Bigl( \fz{Du}{Dt} (x,t) \cdot \nabla \Bigr) \zeta \Bigr] \bigl( x - k \Delta t u(x,t), t-k\Delta t \bigr) = Z(x,t) + O(\Delta t),
	\]
	for the last equality.
\end{proof}
\begin{proof}[Proof of Theorem~\ref{thm:sec_order}]
	In the proof, we often employ simple notations, $L(\cdot,\cdot) = L(x,t;\,\cdot\,, \cdot\,)$ and $X = X(x,t; \,\cdot\,)$, if there is no confusion, since $(x,t)$ is considered as a fixed position in space and time.
	From the Adams--Bashforth method, i.e., for a smooth function~$g$ defined in $\R$, $g^\prime (s) = \fz{1}{2\Delta t} [ 3g(s)-4g(s-\Delta t)+g(s-2\Delta t) ] + O(\Delta t^2)$, we have
	\begin{align}
	\uctd{\zeta} (x,t) 
	& = (\mcL_u\zeta) (x,t)
	= (\mcL_u\zeta)\bigl( X(s), s \bigr)_{|s=t} 
	= L(t,s) \prz{}{s}\Bigl[ L(s,t) \, \zeta \bigl( X(s), s \bigr) L(s,t)^{ \top } \Bigr] L(t,s)^\top {}_{{\dps |}s=t} \notag\\
	& = L(t,s) \fz{1}{2\Delta t}\Bigl[
	3 L(s,t) \, \zeta \bigl( X(s), s \bigr) L(s,t)^{ \top } 
	-4 L(s-\Delta t,t) \, \zeta \bigl( X(s-\Delta t), s-\Delta t \bigr) L(s-\Delta t,t)^{ \top } \notag\\
	& \quad + L(s-2\Delta t,t) \, \zeta \bigl( X(s-2\Delta t), s-2\Delta t \bigr) L(s-2\Delta t,t)^{ \top } \Bigr] L(t,s)^{ \top } {}_{{\dps |}s=t} + O(\Delta t^2) \notag\\
	& = \fz{1}{2\Delta t}\Bigl[
	3 \zeta \bigl( X(x,t; s), s \bigr) 
	-4 L(t,s) L(s-\Delta t,t) \, \zeta \bigl( X(s-\Delta t), s-\Delta t \bigr) L(s-\Delta t,t)^{ \top } L(t,s)^{ \top }  \notag\\
	& \quad + L(t,s)L(s-2\Delta t,t) \, \zeta \bigl( X(s-2\Delta t), s-2\Delta t \bigr) L(s-2\Delta t,t)^{ \top } L(t,s)^{ \top }
	\Bigr] {}_{{\dps |}s=t} + O(\Delta t^2) \quad \mbox{(by \eqref{eqns:L1})} \notag\\
	& = \fz{1}{2\Delta t}\Bigl[
	3 \zeta (x,t) 
	-4 L(t-\Delta t,t) \, \zeta \bigl( X(t-\Delta t), t-\Delta t \bigr) L(t-\Delta t,t)^{ \top }  \notag\\
	& \quad + L(t-2\Delta t,t) \, \zeta \bigl( X(t-2\Delta t), t-2\Delta t \bigr) L(t-2\Delta t,t)^{ \top }
	\Bigr] + O(\Delta t^2) \quad \mbox{(by \eqref{eq:ode2} and \eqref{eqns:L1})} \notag\\
	& = \fz{1}{2\Delta t}\biggl[
	3 \zeta (x,t) -4 \Bigl[ L_1 + \fz{\Delta t^2}{2} U \Bigr] (x,t) \, \zeta \bigl( X(t-\Delta t), t-\Delta t \bigr) \Bigl[ L_1 + \fz{\Delta t^2}{2} U \Bigr]^{ \top }(x,t)  \notag\\
	& \quad + \bigl[ \tilde{L}_1 + 2 \Delta t^2 U \bigr](x,t) \, \zeta \bigl( X(t-2\Delta t), t-2\Delta t \bigr) \bigl[ \tilde{L}_1 + 2 \Delta t^2 U \bigr]^{ \top } (x,t)
	\biggr] + O(\Delta t^2) \notag\\
	& \qquad\qquad\qquad\qquad\qquad\qquad\qquad\qquad\qquad\qquad \mbox{(by Lem.~\ref{lem:L_with_k} with definitions of $L_1$ and $\tilde{L}_1$)} \notag\\
	& = \fz{1}{2\Delta t}\Bigl[
	3 \zeta (x,t) -4 L_1 (x,t) \, \zeta \bigl( X(t-\Delta t), t-\Delta t \bigr) L_1^\top (x,t)  \notag\\
	& \quad + \tilde{L}_1 (x,t) \, \zeta \bigl( X(t-2\Delta t), t-2\Delta t \bigr) \tilde{L}_1^\top (x,t) \notag\\
	& \quad - 2\Delta t^2 \bigl[ \zeta \bigl( X(t-\Delta t), t-\Delta t \bigr) - \zeta \bigl( X(t-2\Delta t), t-2\Delta t \bigr) \bigr] U^{ \top }(x,t) \notag\\
	& \quad - 2\Delta t^2 U(x,t) \bigl[ \zeta \bigl( X(t-\Delta t), t-\Delta t \bigr) - \zeta \bigl( X(t-2\Delta t), t-2\Delta t \bigr) \bigr]
	\Bigr] + O(\Delta t^2) \notag\\
	& = \fz{1}{2\Delta t}\Bigl[
	3 \zeta (x,t) - 4 L_1 (x,t) \, \zeta \bigl( X(t-\Delta t), t-\Delta t \bigr) L_1^\top(x,t) \notag\\
	& \quad + \tilde{L}_1(x,t) \, \zeta \bigl( X(t-2\Delta t), t-2\Delta t \bigr) \tilde{L}_1^\top (x,t)
	\Bigr] + O(\Delta t^2),
	\label{eq:proof1}
	\end{align}
	where the relation,
	\[
	\zeta \bigl( X(t-\Delta t), t-\Delta t \bigr) - \zeta \bigl( X(t-2\Delta t), t-2\Delta t \bigr) = O(\Delta t),
	\]
	has been employed for the last equality.
	Combining~Lemma~\ref{lem:phi_g_with_k}-(ii) with~\eqref{eq:proof1} and recalling $x-\Delta t u(x,t) = X_1(x,t)$ and $x-2\Delta t u(x,t) = \tilde{X}_1(x,t)$, we obtain
	\begin{align*}
	\uctd{\zeta} (x,t) = \fz{1}{2\Delta t}\Bigl[
	3 \zeta (x,t) 
	-4 L_1 (x,t) \zeta \bigl( X_1(x,t), t-\Delta t \bigr) L_1^\top(x,t) \\
	+ \tilde{L}_1 (x,t) \zeta \bigl( \tilde{X}_1(x,t), t-2\Delta t \bigr) \tilde{L}_1^\top (x,t)
	\Bigr] + O(\Delta t^2),
	\end{align*}
	which completes the proof.
\end{proof}
Substituting $t^n$ into $t$ in~\eqref{eq:sec_order}, the discrete form of second-order in time for the upper-convected time derivative is given as follows.
\begin{Cor}\label{cor:sec_order}
Under the same assumptions of Theorem~\ref{thm:sec_order}, we have
\begin{align}
	\uctd{\zeta} (x, t^n) 
	= \fz{1}{2\Delta t}\Bigl[
	3\zeta^n(x) -4 L_1^n (x) \bigl( \zeta^{n-1}\circ X_1^n \bigr)(x) L_1^n (x)^\top 
	+ \tilde{L}_1^n (x) \bigl( \zeta^{n-2}\circ \tilde{X}_1^n \bigr)(x) \tilde{L}_1^n (x)^\top \Bigr] + O(\Delta t^2) 
	\label{eq:sec_order_n}
	\end{align}
for $n = 2, \ldots, N_T$.
\end{Cor}
\begin{Rmk}
Although the approximation~\eqref{eq:sec_order_n} of $\uctd{\zeta}(x,t^n)$ of second-order in time is combined with the finite difference method in this paper below, one can combine it with other methods, e.g., the finite element method and the finite volume method.
\end{Rmk}
%
%
%
%
%
%
%
\subsection{Full discretizations of the upper-convected time derivative}
%
%
%
%
Suppose that $\zeta \in C([0,T];C(\bar\Omega;\R^{d\times d}_\sym))$ and $\zeta_h = \{\zeta_h^n\}_{n=0}^{N_T}\subset V_h$ are given.
For $n \in \{1, \ldots, N_T\}$ and $p\in\{1, 2\}$, let $\mathcal{A}^n \zeta: \bar\Omega \to \R^{d\times d}_\sym$ and $\mathcal{A}_h^{n,(p)} \zeta_h: \bar\Omega_h \to \R^{d\times d}_\sym$ be functions defined by
\begin{align}
	[\mathcal{A}^n \zeta] (x)
	& \defeq 
	\left\{
	\begin{aligned}
	\fz{1}{2\Delta t}\Bigl[
	3\zeta^n(x) 
	-4 L_1^n (x) \bigl( \zeta^{n-1}\circ X_1^n \bigr)(x) L_1^n (x)^\top \qquad \\
	+ \tilde{L}_1^n (x) \bigl( \zeta^{n-2}\circ \tilde{X}_1^n \bigr)(x) \tilde{L}_1^n (x)^\top \Bigr] \quad (n \ge 2), \\
	\fz{1}{\Delta t} \Bigl[
	\zeta^1(x) 
	- L_1^1 (x) \bigl( \zeta^0 \circ X_1^1 \bigr) (x) L_1^1 (x)^\top \Bigr] \quad (n = 1),
	\end{aligned}
	\right. \notag \\
	\label{def:An_h}
	[\mathcal{A}_h^{n,(p)} \zeta_h] (x) 
	& \defeq \left\{
	\begin{aligned}
	\fz{1}{2\Delta t}  \Bigl[
	3\zeta_h^n (x)
	-4 L_1^n (x) \bigl[ \bigl( \Pi_h^{(p)}\zeta_h^{n-1}\bigr) \circ X_1^n \bigr] (x) L_1^n (x)^\top \qquad \\
	+ \tilde{L}_1^n (x) \bigl[ \bigl( \Pi_h^{(p)}\zeta_h^{n-2} \bigr) \circ \tilde{X}_1^n \bigr] (x) \tilde{L}_1^n (x)^\top
	\Bigr] \quad (n \ge 2), \\
	\fz{1}{\Delta t} \Bigl[
	\zeta_h^1 (x)
	- L_1^1 (x) \bigl[ \bigl( \Pi_h^{(p)}\zeta_h^0 \bigr) \circ X_1^1 \bigr] (x) L_1^1 (x)^\top \Bigr] \quad (n = 1),
	\end{aligned}
	\right.
\end{align}
respectively.
Using the notation $\mathcal{A}^n \zeta$, we can write Eq.~\eqref{eq:sec_order_n} as, for $n=\{2,\ldots,N_T\}$,
\[
	\uctd{\zeta} (x, t^n) = [\mathcal{A}^n \zeta] (x) + O(\Delta t^2).
\]
\par
Now, we present a theorem on the truncation error of our finite difference approximations of the upper-convected time derivative, where the function~$\mathcal{A}_h^{n,(p)} \zeta: \bar\Omega_h \to \R^{d\times d}_\sym$ to be used in the theorem has meaning since $\zeta \in C([0,T];C(\bar\Omega;\R^{d\times d}_\sym))$ can be considered as a series of functions in $V_h$, i.e., $\zeta = \{\zeta^n\}_{n=0}^{N_T} \subset V_h$.
\begin{Thm}\label{prop:sec_order_p_order}
	Suppose that Hypotheses~\ref{hyp:u},~\ref{hyp:dt} and~\ref{hyp:mesh} hold true.
	Let $\zeta: \bar\Omega\times [0,T] \to \R^{d\times d}$ be a sufficiently smooth function.
	Then, we have
	\begin{equation}
	\label{eq:sec_order_p_order}
	\uctd{\zeta} (x, t^n) = [\mathcal{A}_h^{n,(p)} \zeta] (x) + O(\Delta t^2 + h^p)
	\end{equation}
	for $x\in\bar\Omega_h$, $n \in \{2, \ldots, N_T\}$ and $p=1, 2$.
\end{Thm}
\begin{proof}
Since for $x\in\bar\Omega_h$ we have
\begin{align}
	\uctd{\zeta} (x, t^n)
	& = [\mathcal{A}^n \zeta] (x) + O(\Delta t^2)  \notag\\
	& = [\mathcal{A}_h^{n,(p)} \zeta] (x) - \Bigl( [\mathcal{A}_h^{n,(p)} \zeta] (x) - [\mathcal{A}^n \zeta] (x) \Bigr) + O(\Delta t^2) \notag\\
	& = [\mathcal{A}_h^{n,(p)} \zeta] (x) + \fz{2}{\Delta t} L_1^n (x) \bigl[ \bigl( \Pi_h^{(p)}\zeta^{n-1} - \zeta^{n-1}\bigr) \circ X_1^n \bigr] (x) L_1^n (x)^\top \notag \\
	& \quad - \fz{1}{2\Delta t} \tilde{L}_1^n (x) \bigl[ \bigl( \Pi_h^{(p)}\zeta^{n-2} - \zeta^{n-2} \bigr) \circ \tilde{X}_1^n \bigr] (x) \tilde{L}_1^n (x)^\top + O(\Delta t^2)
	\label{eq:prop2_proof0}
\end{align}
from Corollary~\ref{cor:sec_order}, it is enough for the proof to show the following estimates,
\begin{subequations}
\label{eq:prop2_proof1}
\begin{align}
\fz{2}{\Delta t}\bigl[ \bigl( \Pi_h^{(p)}\zeta^{n-1} - \zeta^{n-1}\bigr) \circ X_1^n \bigr] (x) & = O(h^p), \\
\fz{1}{2\Delta t}\bigl[ \bigl( \Pi_h^{(p)}\zeta^{n-2} - \zeta^{n-2} \bigr) \circ \tilde{X}_1^n \bigr] (x)
& = O(h^p),
\end{align}
\end{subequations}
where simple estimates~\eqref{eqns:prop2_proof1_easy} are easily obtained as shown in Remark~\ref{rmk:prop2_proof1_easy_interpolation} later and the key issue is to eliminate the negative order in $\Delta t$ from~\eqref{eqns:prop2_proof1_easy} and get~\eqref{eq:prop2_proof1}.
We prove the former equality of~\eqref{eq:prop2_proof1} for $d=2$ only, as the equality for $d=1$ is simpler and the latter one is proved similarly.
Let $x = x_{i,j} \in \bar\Omega_h$ and $y^n \defeq X_1^n(x) = x - u^n(x)\Delta t$.
To simplify notations, we omit superscripts ${}^{n-1}$ and ${}^n$ from $\zeta^{n-1}$ and $y^n$ in the rest of proof, respectively, if there is no confusion.
\par
Let us start with $p=1$.
From Hypotheses~\ref{hyp:u} and~\ref{hyp:dt}, we have $y \in\bar\Omega$ and there exists a pair of indexes~$(i_0, j_0)$ such that $y \in K^{(1)}_{i_0+1/2,j_0+1/2} \, (= [i_0h_1, (i_0+1)h_1] \times [j_0h_2, (j_0+1)h_2])$.
Let $\Lambda^{(1)}(y)$ be a set of pairs of indexes of lattice points near $y$ defined by $\Lambda^{(1)}(y) \defeq \{ (i_0,j_0), (i_0+1,j_0), (i_0,j_0+1), (i_0+1,j_0+1) \}$.
Let $a=(a_1,a_2)^{ \top } \defeq y-x_{i_0,j_0} = ((i-i_0)h_1-u^n_1(x_{i,j})\Delta t, (j-j_0)h_2-u^n_2(x_{i,j})\Delta t)$ and $\tilde{a}=(\tilde{a}_1,\tilde{a}_2)^{ \top } \defeq x_{i_0+1,j_0+1}-y$.
Without loss of generality, we can assume that $u^n_k(x_{i,j}) \ge 0~(k=1,2)$, $i_0 < i$, $j_0 < j$ and $a_k, \tilde{a}_k \ge 0~(k=1,2)$, cf. Fig.~\ref{fig:a1a2}.
\begin{figure}[!htbp]
\centering
\includegraphics[width=0.6\linewidth]{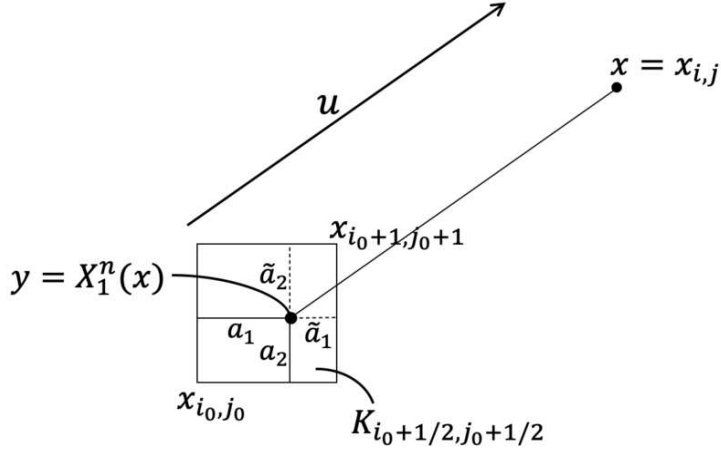}
\caption{Notations in the proof of Theorem~\ref{prop:sec_order_p_order}}
\label{fig:a1a2}
\end{figure}
Then, we have
\begin{align}
& \Bigl[ \bigl( \Pi_h^{(1)}\zeta - \zeta \bigr) \circ X_1^n \Bigr] (x)
= (\Pi_h^{(1)}\zeta) (y) - \zeta (y) \notag\\
& = \sum_{(k,l) \in \Lambda^{(1)}(y)} \bigl[ \zeta (x_{k,l}) - \zeta (y) \bigr] \varphi^{(1)}_{k,l} (y) 
\qquad \mbox{(by~$\sum_{(k,l)\in\Lambda^{(1)}(y)} \varphi^{(1)}_{k,l} (y) = 1$)} \notag\\
& = \sum_{(k,l) \in \Lambda^{(1)}(y)} \Bigl[ \zeta \bigl( y + s(x_{k,l}-y) \bigr) \Bigr]_{s=0}^1 \, \varphi^{(1)}_{k,l} (y) \notag\\
& = \sum_{(k,l) \in \Lambda^{(1)}(y)} \int_0^1 \Bigl( [ (x_{k,l}-y)\cdot\nabla ] \zeta \Bigr) \bigl( y + s_1(x_{k,l}-y) \bigr) ds_1 \, \varphi^{(1)}_{k,l} (y) \notag\\
\label{eq:prop2_proof2}
& = \sum_{(k,l) \in \Lambda^{(1)}(y)} \int_0^1ds_1\int_0^{s_1} \Bigl( [ (x_{k,l}-y)\cdot\nabla ]^2 \zeta \Bigr) \bigl( y + s_2(x_{k,l}-y) \bigr) ds_2 \, \varphi^{(1)}_{k,l} (y) \\
& \qquad\qquad\qquad\qquad\qquad \mbox{(by~$\sum_{(k,l) \in \Lambda^{(1)}(y)} ( [ (x_{k,l}-y)\cdot\nabla ] \zeta ) ( y ) \varphi^{(1)}_{k,l} (y) = 0$),} \notag
\end{align}
and, for $(k,l)=(i_0,j_0)$,
\begin{align}
& \biggl| \int_0^1ds_1\int_0^{s_1} \Bigl( [ (x_{k,l}-y)\cdot\nabla ]^2 \zeta \Bigr) \bigl( y + s_2(x_{k,l}-y) \bigr) ds_2 \, \varphi^{(1)}_{k,l} (y) \biggr| 
= \biggl| \int_0^1ds_1\int_0^{s_1} \bigl( [ a \cdot\nabla ]^2 \zeta \bigr) ( y - s_2 a ) ds_2 \, \fz{\tilde{a}_1\tilde{a}_2}{h_1h_2} \biggr| \notag\\
& \le c_1 (a_1+a_2)^2 \|\zeta^{n-1}\|_{C^2(K^{(1)}_{i_0+1/2,j_0+1/2};\R^{d\times d}_\sym)} \fz{\tilde{a}_1\tilde{a}_2}{h_1h_2} 
\le c_1^\prime (a_1\tilde{a}_1+a_2\tilde{a}_2) \|\zeta^{n-1}\|_{C^2(\bar\Omega;\R^{d\times d}_\sym)} \notag\\ 
& \qquad\qquad\qquad\qquad\qquad\qquad\qquad\qquad\qquad\qquad\qquad\qquad
\mbox{(by $a_k, \tilde{a}_k \le h_k,\ k=1, 2$, and Hyp.~\ref{hyp:mesh})}
\label{eq:prop2_proof3}
\end{align}
for positive constants~$c_1$ and $c_1^\prime$ independent of $h$ and $\Delta t$.
\par
We evaluate $a_1\tilde{a}_1$.
Let $U^\infty \defeq \|u\|_{C([0,T];C(\bar\Omega;\R^d))} = \max\{ |u_k(x,t)|; x \in \bar\Omega, t \in [0,T], k=1,2\}$.
From $y_1 = [x_{i,j} - u^n(x_{i,j})\Delta t ]_1 \in [i_0h_1, (i_0+1)h_1]$, it holds that
\begin{displaymath}
 (i-i_0-1)h_1 \le u^n_1(x_{i,j})\Delta t \le (i-i_0)h_1.
\end{displaymath}
In the case of $i-i_0-1 \in \N$, from $h_1 \le \fz{u^n_1(x_{i,j})\Delta t}{i-i_0-1} \le U^\infty\Delta t$, we have $a_1\tilde{a}_1 \le h_1^2 \le h_1 U^\infty \Delta t$.
In the case of $i-i_0-1 = 0$, from $a_1 \le h_1$ and $\tilde{a}_1 = u^n_1(x_{i,j})\Delta t \le U^\infty\Delta t$, we have $a_1\tilde{a}_1 \le h_1 U^\infty \Delta t$.
Hence, it holds that, for any case,
\[
a_1\tilde{a}_1 \le h_1 U^\infty \Delta t.
\]
Since it holds that $a_2\tilde{a}_2 \le h_2 U^\infty \Delta t$ similarly, we obtain
\begin{equation}
\label{eq:prop2_proof4}
 a_1\tilde{a}_1 + a_2\tilde{a}_2 \le 2 h U^\infty \Delta t,
\end{equation}
where this estimate holds also for $(k,l)=(i_0+1,j_0), (i_0,j_0+1), (i_0+1,j_0+1)$ similarly.
Combining~\eqref{eq:prop2_proof3} and~\eqref{eq:prop2_proof4} with~\eqref{eq:prop2_proof2}, we have, for a positive constant~$c_2$ independent of $h$ and~$\Delta t$,
\[
 \fz{2}{\Delta t}\Bigl[ \bigl( \Pi_h^{(1)}\zeta - \zeta \bigr) \circ X_1^n \Bigr] (x) 
  \le c_2 U^\infty h  \|\zeta\|_{C([0,T];C^2(\bar\Omega;\R^{d\times d}_\sym))} = O(h),
\]
which implies the former equality in~\eqref{eq:prop2_proof1} with $p=1$, and the latter is obtained similarly.
Thus, we get~\eqref{eq:sec_order_p_order} with $p=1$.
\par
In the case of $p=2$, the result, i.e., \eqref{eq:sec_order_p_order} with $p=2$, are obtained similarly by taking into account the next identity,
\begin{align*}
\Bigl[ \bigl( \Pi_h^{(2)}\zeta - \zeta \bigr) \circ X_1^n \Bigr] (x) 
= \sum_{(k,l) \in \Lambda^{(2)}(y)} \int_0^1ds_1 \int_0^{s_1}ds_2 \int_0^{s_2} \Bigl( [ (x_{k,l}-y)\cdot\nabla ]^3 \zeta \Bigr) \bigl( y + s_3(x_{k,l}-y) \bigr) ds_3 \, \varphi^{(2)}_{k,l} (y),
\end{align*}
where $\Lambda^{(2)}(y) \defeq \{ (2i_{\ast}+p,2j_{\ast}+q);\ p, q = 0, 1, 2\}$ for $i_{\ast}\in \{ 0, \ldots, M_1\}$ and $j_{\ast} \in \{ 0, \ldots, M_2\}$ satisfying $y\in [2i_{\ast}h_1, 2(i_{\ast}+1)h_1] \times [2j_{\ast}h_2, 2(j_{\ast}+1)h_2]$.
\end{proof}
\begin{Rmk}
\label{rmk:prop2_proof1_easy_interpolation}
It is obvious that
\begin{align}
\uctd{\zeta} (x, t^n) = [\mathcal{A}_h^{n,(p)} \zeta] (x) + O\Bigl(\Delta t^2 + \fz{h^{p+1}}{\Delta t}\Bigr)
\label{eq:dt2_h_to_p_plus_1_over_dt}
\end{align}
for $x\in\bar\Omega_h$, $n \in \{2, \ldots, N_T\}$ and $p=1, 2$, since $\Pi_h^{(p)}\zeta$ has an accuracy of $O(h^{p+1})$.
In fact, from the approximation property of $\Pi_h^{(p)}\zeta$, we have
\begin{subequations}
\label{eqns:prop2_proof1_easy}
\begin{align}
\fz{2}{\Delta t}\bigl[ \bigl( \Pi_h^{(p)}\zeta^{n-1} - \zeta^{n-1}\bigr) \circ X_1^n \bigr] (x) & = O \Bigl( \fz{h^{p+1}}{\Delta t} \Bigr), \\
\fz{1}{2\Delta t}\bigl[ \bigl( \Pi_h^{(p)}\zeta^{n-2} - \zeta^{n-2} \bigr) \circ \tilde{X}_1^n \bigr] (x)
& = O \Bigl( \fz{h^{p+1}}{\Delta t} \Bigr),
\end{align}
\end{subequations}
and the relation~\eqref{eq:dt2_h_to_p_plus_1_over_dt} is obtained by combining~\eqref{eqns:prop2_proof1_easy} with~\eqref{eq:prop2_proof0}.
Theorem~\ref{prop:sec_order_p_order} eliminates the negative order in $\Delta t$ from~\eqref{eq:dt2_h_to_p_plus_1_over_dt} and ensures that we can take small~$\Delta t$ even for a fixed mesh size from a view point of accuracy.
\end{Rmk}
%
%
%
%
%
%
%
%
\section{Numerical schemes}
\label{nm}
%
%
%
%
In this section, we present finite difference schemes of second-order in time and of first- and second-order in space for problem~\eqref{eqns:prob} by using the ideas of discretizations given in Section~\ref{sec:FD_discretizations}.
\par
Suppose that $u\in C^0([0,T];C^1(\bar\Omega;\R^d))$ and $\zeta^0\in C^0(\bar\Omega;\R^{d\times d}_\sym)$ are given, and that Hypotheses~\ref{hyp:u}, \ref{hyp:dt} and~\ref{hyp:mesh} hold true.
Our schemes are written in a unified form for $d=1, 2(, 3)$ and $p=1, 2$; find $\{\zeta_h^n \in V_h;\ n=1, \ldots, N_T\}$ such that
\begin{subequations}\label{scheme}
	\begin{align}
	\label{scheme_An_h:general_step}
	&&&& [\mathcal{A}_h^{n,(p)} \zeta_h] (x) & = F^n(x), & x &\in \bar\Omega_h, \ n\ge 1, &&&& \\
	\label{scheme_An_h:initial_step}
	&&&& \zeta_h^0(x) & = \zeta^0(x), & x &\in \bar\Omega_h, &&&&
	\end{align}
\end{subequations}
which are equivalent to
\begin{subequations}\label{scheme_detail}
	\begin{align}
	\fz{1}{2\Delta t}  \Bigl[
	3\zeta_h^n(x)
	-4 L_1^n (x)  \bigl[ \bigl( \Pi_h^{(p)}\zeta_h^{n-1}\bigr) \circ X_1^n \bigr] (x) & L_1^n (x)^\top \notag\\
	+ \tilde{L}_1^n (x) \bigl[ \bigl( \Pi_h^{(p)}\zeta_h^{n-2} \bigr) \circ \tilde{X}_1^n \bigr] (x) \tilde{L}_1^n (x)^\top
	\Bigr] & = F^n(x), 
	& x & \in\bar\Omega_h, \ n\ge 2, 
	\label{scheme:general_step} \\
	\label{scheme:first_step}
	\fz{1}{\Delta t} \Bigl[
	\zeta_h^1(x)
	- L_1^1 (x)  \bigl[ \bigl( \Pi_h^{(p)}\zeta_h^0 \bigr) \circ X_1^1 \bigr] (x) L_1^1 (x)^\top \Bigr] 
	& = F^1(x), & x & \in\bar\Omega_h, \\
	\label{scheme:initial_step}
	\zeta_h^0(x) & = \zeta^0(x), & x & \in\bar\Omega_h.
	\end{align}
\end{subequations}
The unified scheme~\eqref{scheme} (equivalent to~\eqref{scheme_detail}) includes four schemes, i.e., $p=1$ and $2$ correspond to schemes of first- and second-order in space, respectively, and the spatial dimension~$d~(= 1, 2)$ is implicitly dealt in the symbols~$\bar\Omega_h$ and~$V_h$.
An approximate initial value~$\zeta_h^0 \in V_h$ is given by~\eqref{scheme:initial_step}.
We find $\zeta_h^1\in V_h$ from~\eqref{scheme:first_step} and  $\zeta_h^n\in V_h$ for $n \ge 2$ from~\eqref{scheme:general_step}.
Here, we additionally provide practical form of~\eqref{scheme}:
\begin{subequations}\label{scheme_explicit}
	\begin{align}
	\zeta_h^n(x) & =
	\fz{4}{3} L_1^n (x) \bigl[ \bigl( \Pi_h^{(p)}\zeta_h^{n-1} \bigr) \circ X_1^n \bigr](x) L_1^n (x)^\top \notag\\
	& \quad 
	- \fz{1}{3} \tilde{L}_1^n (x) \bigl[ \bigl( \Pi_h^{(p)}\zeta_h^{n-2} \bigr) \circ \tilde{X}_1^n \bigr] (x) \tilde{L}_1^n (x)^\top
	+ \fz{2\Delta t}{3} F^n(x), & x & \in \bar\Omega_h, \  n \ge 2, \label{scheme_explicit:general_step} \\
	\label{scheme_explicit:first_step}
	\zeta_h^1(x) & =
	L_1^1 (x) \bigl[ \bigl( \Pi_h^{(p)}\zeta_h^0 \bigr) \circ X_1^1 \bigr] (x) L_1^1 (x)^\top + \Delta t F^1(x), &
	x & \in \bar\Omega_h, \\
	\label{scheme_explicit:initial_step}
	\zeta^0_h(x) & = \zeta^0(x), &
	x & \in \bar\Omega_h,
	\end{align}
\end{subequations}
which imply that scheme~\eqref{scheme} is explicit.
\begin{Rmk}
From Hypotheses~\ref{hyp:u} and~\ref{hyp:dt} and Remark~\ref{rmk:upwind_points}, we have $\Gamma_{\rm in} = \emptyset$ and $X_1(\Omega,t)=\tilde{X}_1(\Omega,t)=\Omega~(t\in [0, T])$, i.e., all of upwind points are in $\bar\Omega$.
Hence, the functions~$(\Pi_h^{(p)}\zeta_h^{n-1}) \circ X_1^n~(n \ge 1)$ and $(\Pi_h^{(p)}\zeta_h^{n-2}) \circ \tilde{X}_1^n~(n \ge 2)$ are well defined in $\bar\Omega$ for $p=1,2$.
\end{Rmk}
\begin{Rmk}
In scheme~\eqref{scheme}, we employ the backward {\rm Euler} method~\eqref{scheme:first_step} of first-order in time once to find $\zeta_h^1$ needed in~\eqref{scheme:general_step} with $n=2$. It is expected that there is no influence on the second-order convergence in time, cf. \cite{Notsu2016}.
\end{Rmk}
\begin{Rmk}\label{rmk:symmetry}
Suppose that Hypotheses~\ref{hyp:u} and~\ref{hyp:dt} hold true.
Then, under $F \in C(\bar\Omega\times [0,T]; \R^{d\times d}_\sym)$ and $\zeta^0 \in C(\bar\Omega; \R^{d\times d}_\sym)$, the scheme~\eqref{scheme} preserves the symmetry, i.e., $\zeta_h^n(x)^{ \top } = \zeta_h^n(x)~(x\in\bar\Omega_h,\ n=0,\ldots,N_T)$ from the following.
For $d=1$, it is obvious, and let us consider $d = 2(, 3)$.
$\zeta_h^0(x)~(x\in\bar\Omega_h)$ is symmetric from the symmetry of~$\zeta^0$.
We show the symmetry of $\zeta^1_h(x)~(x\in\bar\Omega_h)$.
Noting~\eqref{scheme_explicit:first_step} and letting $A(x) = L_1^1(x)$, $B(x) = \bigl[ \bigl( \Pi_h^{(p)}\zeta_h^0 \bigr) \circ X_1^1 \bigr] (x)$, and $C(x) = \Delta t F^1(x)$, we have
	\begin{align*}
	\zeta^1_h(x)^{ \top }
	& = \bigl[ A(x) B(x) A(x)^{ \top } + C(x) \bigr]^{ \top } 
	= A(x) B(x)^{ \top } A(x)^{ \top } + C(x)^{ \top } \\
	& = A(x) B(x) A(x)^{ \top } + C(x)
	= \zeta^1_h(x),
	\end{align*}
which implies symmetry of~$\zeta^1_h(x)$ for~$x\in\bar\Omega_h$, where we have used the fact that $B(x)$ and $C(x)$ are symmetric for the second equality from the last.
For~$n\ge 2$, the symmetry of~$\zeta_h^n(x)$ is obtained similarly from~\eqref{scheme:general_step}.
\end{Rmk}
%
%
%
%
%
\subsection{Schemes in one-dimensional space~$(d=1)$}
%
%
%
%
In this subsection, we rewrite the finite difference scheme~\eqref{scheme} in a unified form for $d=1$ and $p=1, 2$.
We introduce simplified notations, $\zeta^n_i \defeq \zeta_h^n(x_i)$, $u^n_i \defeq u^n(x_i)$, $\nabla u^n_i \defeq (\nabla u^n)(x_i)$, $F^n_i \defeq F(x_{i},t^{n})$, $\Lambda_\Omega \defeq \{0,\ldots,N\}$, and $\Lambda_T \defeq \{1,\ldots,N_T\}$.
The schemes are to find $\{\zeta^n_i \in \R;\ i \in \Lambda_\Omega,\ n \in \Lambda_T \}$ such that
\begin{subequations}\label{scheme_1d}
	\begin{align}
	\zeta^n_i & = \fz{4}{3} ( 1 + \Delta t \nabla u^n_i )^2 \bigl[ \bigl( \Pi_h^{(p)}\zeta_h^{n-1}\bigr) \circ X_1^n \bigr] (x_i) 
	\qquad\qquad\ \notag\\
	& \quad -\fz{1}{3} ( 1 + 2\Delta t \nabla u^n_i )^2 \bigl[ \bigl( \Pi_h^{(p)}\zeta_h^{n-2} \bigr) \circ \tilde{X}_1^n \bigr] (x_i) + \fz{2\Delta t}{3} F^n_{i}, 
	& i & \in \Lambda_\Omega, \ n \ge 2, 
	\label{scheme_1d:general_step} \\
	\zeta^1_i & = ( 1 + \Delta t \nabla u^1_i )^2 \bigl[ \bigl( \Pi_h^{(p)}\zeta_h^0 \bigr) \circ X_1^1 \bigr] (x_i) + \Delta t F^1_i, & i &\in \Lambda_\Omega, 
	\label{scheme_1d:first_step} \\
	\zeta^0_i & = \zeta^0(x_i), & i & \in \Lambda_\Omega.
	\label{scheme_1d:initial_step}
	\end{align}
\end{subequations}
\par
We give the algorithm as follows:
\smallskip
\par
\textbf{Algorithm~1~($d=1$).}
Set $\bar\Omega_h = \{ x_i \in\bar\Omega;\ i \in \Lambda_\Omega \}$ with $h=a/N$, and $\{ \zeta^0_i;\ i\in\Lambda_\Omega\}$ by~\eqref{scheme_1d:initial_step} to get $\zeta_h^0\in V_h$, where $N$ is an even number and $M = N/2$ for $p=2$.
	\smallskip\\
	Set~$n=1$. \smallskip\\
	\quad
	For each $i \in \Lambda_\Omega$ do:
		\begin{enumerate}
		\item Compute $F^1_i$, $u^1_i$, $\nabla u^1_i$, and $y^1_i \defeq X_1^{1} (x_i) = x_{i} - \Delta t\,u^1_i$.
		\item Compute $Z_i^{1,(p)} \defeq [ ( \Pi_h^{(p)}\zeta_h^0) \circ X_1^1 ] (x_i) = ( \Pi_h^{(p)}\zeta_h^0) (y_i^1)$ according to~\eqref{int1} with $i_0 = \mathcal{I}(y_i^1; 0, a, N)$ for $p=1$, or \eqref{int2} with $k_0 = \mathcal{I}(y_i^1; 0, a, M)$ for $p=2$.
		\item Compute $\zeta^1_i$ by~\eqref{scheme_1d:first_step}, which is equivalent to
		\[
		\zeta^1_i = ( 1 + \Delta t \nabla u^1_i )^2 \, Z_i^{1,(p)} + \Delta t F^1_i.
		\]
		\end{enumerate}
	(Here, computation of $\zeta_h^1 \in V_h$ is completed.) \\
	Set~$n=2$. \smallskip\\
		While $n \le N_T$ do:\\
		\quad
		For each $i\in \Lambda_\Omega$ do:
\begin{enumerate}
	\item Compute $F^n_i$, $u^n_i$, $\nabla u^n_i$, $y^n_i \defeq X_1^{n} (x_i) = x_{i} - \Delta t\,u^n_i$, and $\tilde{y}^n_i \defeq \tilde{X}_1^{n} (x_i) = x_{i} - 2\Delta t\,u^n_i$.
		\item Compute $Z_i^{n,(p)} \defeq [ ( \Pi_h^{(p)}\zeta_h^{n-1}) \circ X_1^n ] (x_i) = ( \Pi_h^{(p)}\zeta_h^{n-1}) (y_i^n)$ according to~\eqref{int1} with $i_0 = \mathcal{I}(y_i^n; 0, a, N)$ for $p=1$, or \eqref{int2} with $k_0 = \mathcal{I}(y_i^n; 0, a, M)$ for $p=2$.
		Similarly, compute $\tilde{Z}_i^{n,(p)} \defeq [ ( \Pi_h^{(p)}\zeta_h^{n-2}) \circ \tilde{X}_1^n ] (x_i) = ( \Pi_h^{(p)}\zeta_h^{n-2}) (\tilde{y}_i^n)$.
		\item Compute $\zeta^n_i$ by~\eqref{scheme_1d:general_step}, which is equivalent to
		\[
		\zeta^n_i = \fz{4}{3}( 1 + \Delta t \nabla u^n_i )^2 \, z_i^{n,(p)} - \fz{1}{3}( 1 + 2\Delta t \nabla u^n_i )^2 \, \tilde{Z}_i^{n,(p)} + \fz{2\Delta t}{3} F^n_i.
		\]
\end{enumerate}		
	\quad
	(Computation of $\zeta_h^n \in V_h$ is completed.) \\
	\quad
	Set $n=n+1$.
%
%
%
%
\subsection{Schemes in two-dimensional space~($d=2$)}\label{nm2}
%
%
%
%
Similarly to the previous subsection, we rewrite the unified finite difference scheme~\eqref{scheme} for $d=2$ and $p=1,2$.
\par
Let us introduce simplified notations, $\zeta^n_{i,j} \defeq \zeta_h^n(x_{i,j})$, $u^n_{i,j} \defeq u^n(x_{i,j})$, $\nabla u^n_{i,j} \defeq (\nabla u^n)(x_{i,j})$, $F^n_{i,j} \defeq F(x_{i,j},t^{n})$, and $\Lambda_\Omega \defeq \{(i,j);\ i = 0,\ldots,N_1,\ j=0,\ldots,N_2\}$.
The schemes are to find $\{\zeta^n_{i,j} \in \R^{2\times2}_\sym;\ (i,j) \in \Lambda_\Omega,\ n \in \Lambda_T \}$ such that
	\begin{subequations}\label{scheme_2d}
		\begin{align}
		\zeta_{i,j}^n & =
		\fz{4}{3} \bigl[I + \Delta t (\nabla u^n_{i,j}) \bigr] \bigl[ \bigl( \Pi_h^{(p)}\zeta_h^{n-1} \bigr) \circ X_1^n \bigr](x_{i,j}) \bigl[I + \Delta t (\nabla u^n_{i,j}) \bigr]^{ \top } \notag\\
		& \quad - \fz{1}{3} \bigl[I + 2\Delta t (\nabla u^n_{i,j}) \bigr] \bigl[ \bigl( \Pi_h^{(p)}\zeta_h^{n-2} \bigr) \circ \tilde{X}_1^n \bigr] (x_{i,j}) \bigl[I + 2\Delta t (\nabla u^n_{i,j}) \bigr]^{ \top } + \fz{2\Delta t}{3} F^n_{i,j}, &
		(i,j) & \in \Lambda_\Omega, \ n\ge 2, 
		\label{scheme_2d:general_step} \\
		\zeta_{i,j}^1 & =
		\bigl[I + \Delta t (\nabla u^1_{i,j}) \bigr] \bigl[ \bigl( \Pi_h^{(p)}\zeta_h^0 \bigr) \circ X_1^1 \bigr] (x_{i,j}) \bigl[I + \Delta t (\nabla u^1_{i,j}) \bigr]^{ \top } + \Delta t F^1_{i,j}, & (i,j) & \in\Lambda_\Omega, \label{scheme_2d:first_step} \\
		\zeta^0_{i,j} & = \zeta^0(x_{i,j}), & (i,j) & \in \Lambda_\Omega.
		\label{scheme_2d:initial_step}
		\end{align}
	\end{subequations}
\par
We give an algorithm of schemes~\eqref{scheme_2d} for $d=2$ and $p=1,2$, while the construction is analogous to Algorithm~1 for $d=1$.
	\smallskip \par
	\textbf{Algorithm~2~($d=2$).}
Set $\bar\Omega_h = \{ x_{i,j} \in\bar\Omega;\ (i,j) \in \Lambda_\Omega \}$ with $h_i=a_i/N_i~(i=1,2)$, and $\{ \zeta^0_{i,j};\ (i,j)\in\Lambda_\Omega\}$ by~\eqref{scheme_2d:initial_step} to get $\zeta_h^0\in V_h$, where $N_i~(i=1,2)$ are even numbers and $M_i = N_i/2~(i=1,2)$ for $p=2$.
	\smallskip\\
	Set~$n=1$. \smallskip\\
	\quad
	For each $(i,j) \in \Lambda_\Omega$ do:
		\begin{enumerate}
		\item Compute $F^1_{i,j}$, $u^1_{i,j}$, $\nabla u^1_{i,j}$, and $y^1_{i,j} \defeq X_1^{1} (x_{i,j}) = x_{i,j} - \Delta t\,u^1_{i,j}$.
		\item Compute $Z_{i,j}^{1,(p)} \defeq [ ( \Pi_h^{(p)}\zeta_h^0) \circ X_1^1 ] (x_{i,j}) = ( \Pi_h^{(p)}\zeta_h^0) (y_{i,j}^1)$ according to~\eqref{int3} with $i_0 = \mathcal{I}((y_{i,j}^1)_1;\ 0, a_1, N_1)$ and $j_0 = \mathcal{I}((y_{i,j}^1)_2;\ 0, a_2, N_2)$ for $p=1$, or \eqref{int4} with $k_0 = \mathcal{I}((y_{i,j}^1)_1;\ 0, a_1, M_1)$ and $l_0 = \mathcal{I}((y_{i,j}^1)_2;\ 0, a_2, M_2)$ for $p=2$.
		\item Compute $\zeta^1_{i,j}$ by~\eqref{scheme_2d:first_step}, which is equivalent to
		\[
		\zeta_{i,j}^1 =
		\bigl[I + \Delta t (\nabla u^1_{i,j}) \bigr] Z_{i,j}^{1,(p)} \bigl[I + \Delta t (\nabla u^1_{i,j}) \bigr]^{ \top } + \Delta t F^1_{i,j}.
		\]
		\end{enumerate}
	(Here, computation of $\zeta_h^1 \in V_h$ is completed.) \\
	Set~$n=2$. \smallskip\\
		While $n \le N_T$ do:\\
		\quad
		For each $(i,j) \in \Lambda_\Omega$ do:
\begin{enumerate}
	\item Compute $F^n_{i,j}$, $u^n_{i,j}$, $\nabla u^n_{i,j}$, $y^n_{i,j} \defeq X_1^{n} (x_{i,j}) = x_{i,j} - \Delta t\,u^n_{i,j}$, and $\tilde{y}^n_{i,j} \defeq \tilde{X}_1^{n} (x_{i,j}) = x_{i,j} - 2\Delta t\,u^n_{i,j}$.
		\item Compute $Z_{i,j}^{n,(p)} \defeq [ ( \Pi_h^{(p)}\zeta_h^{n-1}) \circ X_1^n ] (x_{i,j}) = ( \Pi_h^{(p)}\zeta_h^{n-1}) (y_{i,j}^n)$ according to~\eqref{int3} with $i_0 = \mathcal{I}((y_{i,j}^n)_1;\ 0, a_1, N_1)$ and $j_0 = \mathcal{I}((y_{i,j}^n)_2;\ 0, a_2, N_2)$ for $p=1$, or \eqref{int4} with $k_0 = \mathcal{I}((y_{i,j}^n)_1; 0, a_1, M_1)$ and $l_0 = \mathcal{I}((y_{i,j}^n)_2; 0, a_2, M_2)$ for $p=2$.
		Similarly, compute $\tilde{Z}_{i,j}^{n,(p)} \defeq [ ( \Pi_h^{(p)}\zeta_h^{n-2}) \circ \tilde{X}_1^n ] (x_{i,j}) = ( \Pi_h^{(p)}\zeta_h^{n-2}) (\tilde{y}_{i,j}^n)$.
		\item Compute $\zeta^n_{i,j}$ by~\eqref{scheme_2d:general_step}, which is equivalent to
		\begin{align*}
		\zeta_{i,j}^n & =
		\fz{4}{3} \bigl[I + \Delta t (\nabla u^n_{i,j}) \bigr] Z_{i,j}^{n,(p)} \bigl[I + \Delta t (\nabla u^n_{i,j}) \bigr]^\top - \fz{1}{3} \bigl[I + 2\Delta t (\nabla u^n_{i,j}) \bigr] \tilde{Z}_{i,j}^{n,(p)} \bigl[I + 2\Delta t (\nabla u^n_{i,j}) \bigr]^{ \top } + \fz{2\Delta t}{3} F^n_{i,j}.
		\end{align*}
\end{enumerate}		
	\quad
	(Computation of $\zeta_h^n \in V_h$ is completed.) \\
	\quad
	Set $n=n+1$.
%
%
%
%
\section{Numerical results}
\label{numerics}
%
%
%
%
In this section, numerical results for problems with manufactured solutions are presented to observe experimental convergence orders of proposed schemes.
In the following, we denote scheme~\eqref{scheme} with $p=1$ and $p=2$ by (S1) and (S2), respectively.
From Theorem~\ref{prop:sec_order_p_order}, the expected orders of convergence are of $O(\Delta t^2+h^p)$ for $p=1,2$.
To see the experimental orders of convergence, the efficient choices of $\Delta t$ for (S1) and (S2) are respectively $\Delta t = c\sqrt{h}$ and $\Delta t = c^\prime h$, for positive constants~$c$ and $c^\prime$.
The choices of $\Delta t$ for (S1) and (S2) lead to an expected order of convergence of $O(\Delta t^2)\ (=O(h^p))$.
In the computations below, as mentioned in Remark~\ref{rmk:upwind_cell_1d}-(iii) and Remark~\ref{rmk:upwind_cell_2d}-(iii), we employ a value of $\zeta_{\rm in}$ at closest lattice point to an upwind point~$X_1^n(x)$ or $\tilde{X}_1^n(x)$ for $x=x_{i}~(d=1)$ or $x_{i,j}~(d=2)$ when the upwind point is outside the domain, where the integer-valued index indicator function~$\mathcal{I}$ given by~\eqref{def:index_func} is used.
\par
For $\psi_h:\bar\Omega_h\to\R$ and $\phi_h = \{\phi_h^n:\bar\Omega_h\to\R;\ n=1,\ldots,N_T \}$, let $\|\cdot\|_{\ell^\infty(\bar\Omega_h)}$ and $\|\cdot\|_{\ell^\infty(\ell^\infty)}$ be norms defined by
\begin{align*}
\| \psi_h \|_{\ell^\infty(\bar\Omega_h)} & = \| \psi_h \|_{\ell^\infty(\bar\Omega_h;\R)} \defeq \max \bigl\{ | \psi_h(x) | ;\ x\in \bar\Omega_h \bigr\}, \\
\| \phi_h \|_{\ell^\infty(\ell^\infty)} & \defeq \max \bigl\{ \| \phi_h^n \|_{\ell^\infty(\bar{\Omega}_h)};\ n = 1,\ldots, N_T \bigr\}.
\end{align*}
Let $E_{ij} = E_{ij}(\Delta t, h),\ i,j=1,\ldots,d,$ be errors between a numerical solution~$\zeta_h = \{ \zeta_h^n \}_{n=1}^{N_T} \subset V_h$ and a corresponding exact solution~$\zeta \in C([0,T];C(\bar\Omega;\R^{d\times d}_\sym))$ defined by
\[
 E_{ij}  = E_{ij}(\Delta t, h) \defeq \bigl\| [\zeta_h]_{ij}-\zeta_{ij} \bigr\|_{\ell^\infty(\ell^\infty)}, \quad i, j = 1, \ldots, d,
\]
and $E_{11}$ is denoted by $E$ simply when $d=1$.
\begin{Rmk}
To solve the problems proposed in this Section, we are assuming a defined source term and a prescribed velocity field. In addition, we need to establish at least one initial condition and a wall condition where we can call the flow inlet. The initial condition $\zeta_{h}^{0}(x)$ is directly derived from the exact solution $\zeta_{\rm exact}(x,0)$. The boundary condition is computed assuming a Dirichlet-type condition, i.e., we use the exact solution $\zeta_{\rm in} = \zeta_{\rm exact}^{n}(x_{0})$ at the first point of the boundary for a positive velocity field (if $u^n(x)<0$ then the inlet of the domain is located on the opposite side, making us consider $\zeta_{\rm exact}^{n}(x_{N})$). Therefore, when we have the case described in Fig.~\ref{fig:BC1}, the interpolated point $X^n_1(x_0)$ at previous time is outside the domain; thus we have imposed the boundary condition $\zeta_{\rm in}^{n}(x)=\zeta_{\rm exact}^{n}(x)$. 

	For the opposite side of the domain as represented by Fig.~\ref{fig:NeumannBC1}, we do not impose any wall conditions, since our method can also be used to update the value of unknown function $\zeta_{\rm in}^{n}(x_N)$ on the outflow wall. In addition, it is also possible to assume a Neumann boundary condition on this wall and then we apply the method until $x_{N-1}$ and update the last point as in an explicit scheme $\zeta_{\rm in}^{n}(x_N) = \zeta_{\rm in}^{n}(x_{N-1})$.

	More details about the implementation of these strategies can be found in Appendix~\ref{A.subsec:PseudoCode}.

		\begin{figure}[htbp]
			\centering
			\begin{tikzpicture}[thick,scale=1.7, every node/.style={scale=1.2}]
			\draw[thick,->] (-1.0,-1.0) -- (1.25,-1.0);
			\draw[thick,->] (-1.0,-1.0) -- (-1.0,1.1);
			\draw[thick] (1.0,-1.05) -- (1.0,-0.95);
			
			\draw (1.5,-0.5) node{$t^{n-1}$};
			\draw (1.5,0.0) node{$t^{n}$};
			\draw (-1.05,-1.25) node{$x_0=0$};
			\draw (0.95,-1.25) node{$x_N=a$};
			
			\draw[dashed] (-1.0,-0.5) -- (1.0,-0.5);
			\draw[dashed] (-1.0,0.0) -- (1.0,0.0);
			\draw[dashed,blue] (-1.5,-1.0) -- (-0.5,1.0);
			\draw (0.3,1.0) node{$characteristic$};
			\filldraw[fill] (-1.0,0.0) circle (0.02cm)
			(-1.0,-0.5) circle (0.02cm)
			(-1.25,-0.5) circle (0.02cm);
			\draw (-0.7,0.2) node{$\zeta^n_{\rm in}$};
			\draw (-0.7,-0.3) node{$\zeta^{n-1}_{\rm in}$};
			\draw (-1.75,-0.4) node{$X^{n}_{1}(x_0)$};
			
			\end{tikzpicture}
			\caption{Sketch of the wall treatment for unknown boundary condition.}
			\label{fig:BC1}
		\end{figure}
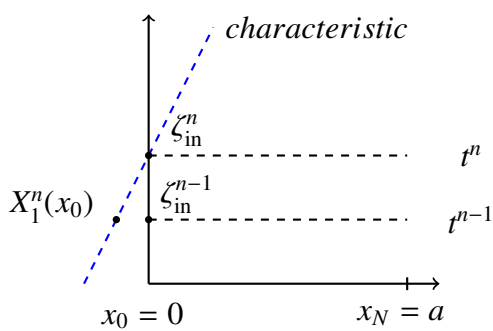

		\begin{figure}[htbp]
			\centering
			\begin{tikzpicture}[thick,scale=1.7, every node/.style={scale=1.7}]
			\draw[thick,->] (-1.0,-1.0) -- (1.25,-1.0);
			\draw[thick,->] (-1.0,-1.1) -- (-1.0,1.1);
			\draw[thick] (1.0,-1.1) -- (1.0,-0.9);
			
			\draw (-1.25,1.0) node{$t$};
			\draw (1.35,-1.25) node{$x$};
			\draw (1.75,0.25) node{$\zeta^n_{i}=\zeta^n_{i-1}$};
			\draw (0.45,0.25) node{$\zeta^n_{i-1}$};
			
			\draw[dashed] (-1.0,-0.5) -- (1.0,-0.5);
			\draw[dashed] (-1.0,0.0) -- (1.0,0.0);
			\draw[dashed,blue] (0.1,-1.0) -- (0.85,0.5);
			\draw[dashed,blue] (0.5,-1.0) -- (1.0,0.0);
			
			\draw (-1.0,-1.2) node{$0$};
			\draw (0.9,-1.2) node{$x_N$};
			%
			\filldraw[fill] (0.35,-0.5) circle (0.02cm)
			(0.6,0.0) circle (0.02cm)
			(1.0,0.0) circle (0.02cm);
			%
			\draw (-1.75,0.1) node{$t^n$};
			\draw (-1.5,-0.4) node{$t^{n-1}$};
			
			\end{tikzpicture}
			\caption{Sketch of the wall treatment for Neumann boundary condition.}
			\label{fig:NeumannBC1}
		\end{figure}
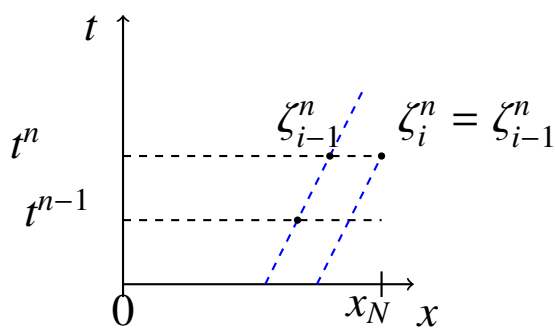

\end{Rmk}
\subsection{Examples in one-dimensional space~$(d=1)$}
%
%
%
%
We consider the next example in one-dimensional space.
\begin{Ex}[$d=1$]\label{ex:1d}
In problem~\eqref{eqns:prob}, let $d=1$, $\Omega = (0,1)$ and $T=1$.
We consider three functions for the velocity:
\[
 (i)~u(x,t)  = t, \qquad
 (ii)~u(x,t)  = x+t,  \qquad
 (iii)~u(x,t)  = \sin(x+t),
\]
which imply $\Gamma_{\rm in} = \{ 0 \}~(t\in (0,T])$.
The functions~$F$, $\zeta_{\rm in}$ and $\zeta^0$ are given so that the exact solution is
\[
 \zeta(x,t) = \sin (x+t) + 2.
\]
\end{Ex}
\par
We solve Example~\ref{ex:1d} by (S1) with $\Delta t = c\sqrt{h}$ for $c = 1/50$ and (S2) with $\Delta t = c^\prime h$ for $c^\prime = 1$, where the mesh is constructed for $h=1/N$ with $N = 10, 20, 40, 80, 160$ and $320$, the constants $c$ and $c^\prime$ are as larger as possible in order to numerically verify the convergence order of the temporal discretizations.
Tables~\ref{table:ex1_s1_dt_sqrth} and \ref{table:ex1_s2_dt_h} show the values of error $E$ and their slopes in~$\Delta t$.
According to the results in the tables,	we can confirm that (S1) and (S2) are of second-order in $\Delta t$ for the three cases of velocity, $(i)$, $(ii)$ and $(iii)$.
These results are consistent with the theoretical results in Theorem~\ref{prop:sec_order_p_order}.
	\begin{table}[!htbp]
		\centering
		\caption{Example~\ref{ex:1d} by {\rm (S1)} with $\Delta t = c \sqrt{h}$~$(c=1/50)$: Values of~$E$ and their slopes in $\Delta t$.}
		\begin{tabular}{rrrcrrcrr}
			\toprule
			& \multicolumn{2}{c}{$(i)$} && \multicolumn{2}{c}{$(ii)$} && \multicolumn{2}{c}{$(iii)$} \\ \cline{2-3}\cline{5-6}\cline{8-9}
			\multicolumn{1}{c}{$N$} & \multicolumn{1}{c}{$E$} & \multicolumn{1}{c}{Slope} && \multicolumn{1}{c}{$E$} & \multicolumn{1}{c}{Slope} && \multicolumn{1}{c}{$E$} & \multicolumn{1}{c}{Slope} \\
			\hline
			$10$ & $1.54 \times 10^{-2}$ & -- && $3.45 \times 10^{-2}$ & -- && $2.11 \times 10^{-2}$ & -- \\
			$20$ & $8.07 \times 10^{-3}$ & $1.86$ && $1.83 \times 10^{-2}$ & $1.87$ && $1.11 \times 10^{-2}$ & $1.86$ \\
			$40$ & $4.15 \times 10^{-3}$ & $1.92$ && $9.38 \times 10^{-3}$ & $1.92$ && $5.69 \times 10^{-3}$ & $1.92$ \\
			$80$ & $2.10 \times 10^{-3}$ & $1.96$ && $4.75 \times 10^{-3}$ & $1.96$ && $2.88 \times 10^{-3}$ & $1.96$ \\
			$160$ & $1.06 \times 10^{-3}$ & $1.98$ && $2.39 \times 10^{-3}$ & $1.98$ && $1.45 \times 10^{-3}$ & $1.98$ \\
			$320$ & $5.31 \times 10^{-4}$ & $1.99$ && $1.13 \times 10^{-3}$ & $2.16$ && $7.27 \times 10^{-4}$ & $1.99$ \\
			\bottomrule
		\end{tabular}
		\label{table:ex1_s1_dt_sqrth}
	\end{table}
	\begin{table}[!htbp]
		\centering
		\caption{Example~\ref{ex:1d} by {\rm (S2)} with $\Delta t = h$~$(c^\prime = 1)$: Values of~$E$ and their slopes in $\Delta t$.}
		\begin{tabular}{rrrcrrcrr}
			\toprule
			& \multicolumn{2}{c}{$(i)$} && \multicolumn{2}{c}{$(ii)$} && \multicolumn{2}{c}{$(iii)$} \\ \cline{2-3}\cline{5-6}\cline{8-9}
			\multicolumn{1}{c}{$N$} & \multicolumn{1}{c}{$E$} & \multicolumn{1}{c}{Slope} && \multicolumn{1}{c}{$E$} & \multicolumn{1}{c}{Slope} && \multicolumn{1}{c}{$E$} & \multicolumn{1}{c}{Slope} \\
			\hline
			$10$ & $4.65 \times 10^{-3}$ & -- && $8.05 \times 10^{-2}$ & -- && $1.65 \times 10^{-2}$ & -- \\
			$20$ & $1.11 \times 10^{-3}$ & $2.07$ && $2.19 \times 10^{-2}$ & $1.88$ && $5.45 \times 10^{-3}$ & $1.60$ \\
			$40$ & $2.68 \times 10^{-4}$ & $2.04$ && $5.63 \times 10^{-3}$ & $1.96$ && $1.53 \times 10^{-3}$ & $1.84$ \\
			$80$ & $6.59 \times 10^{-5}$ & $2.03$ && $1.42 \times 10^{-3}$ & $1.98$ && $4.02 \times 10^{-4}$ & $1.93$ \\
			$160$ & $1.63 \times 10^{-5}$ & $2.01$ && $3.58 \times 10^{-4}$ & $1.99$ && $1.03 \times 10^{-4}$ & $1.97$ \\
			$ 320 $ & $ 4.06 \times 10^{-6}$ & $ 2.01 $ && $8.96 \times 10^{-5}$ & $2.00 $ && $ 2.61 \times 10^{-5}$ & $1.98 $ \\
			\bottomrule
		\end{tabular}
		\label{table:ex1_s2_dt_h}
	\end{table}	

In order to numerically verify that our methodology is stable for small time-steps, we have fixed a coarse mesh $h=1/40$ and the finest mesh $h=1/320$ simulating the reduction of the time-step as $\Delta t(k)=\frac{\sqrt{h}/50}{2^{k}}$ for the first-order scheme and as $\Delta t(k)=\frac{h}{2^{k}}$ for the second-order method. Results for (S1) in Table \ref{tableX1} while in Table \ref{tableX2} we have described the results for (S2). 
	
According to these tables, we can confirm that our methodologies, first- and second-order spatial discretization schemes, are unconditionally stable since the errors are decreasing as $\Delta t$ is reduced. It is important to highlight that error for the smallest time-step in Table \ref{tableX2} for $h=1/40$ is approximately two order smaller than the error of the largest time-step, confirming the good stability property of the second-order scheme.

	\begin{table}[ht!]
	\centering
	\caption{Example~\ref{ex:1d} by {\rm (S1)}: reducing the time-step as $\Delta t(k)=\frac{\sqrt{h}/50}{2^{k}}$.}
	\label{tableX1}
	\begin{tabular}{rrr}
		\hline 
		\multicolumn{3}{c}{$h = 0.025$} \\ \hline
		\multicolumn{1}{c}{$k$} & \multicolumn{1}{c}{$\Delta t$} & \multicolumn{1}{c}{$Error$}  \\
		\hline
		$0$ & $3.16 \times 10^{-3}$ & $1.18375 \times 10^{-2}$   \\
		$1$ & $1.58 \times 10^{-3}$ & $1.08387 \times 10^{-2}$   \\
		$2$ & $7.91 \times 10^{-4}$ & $1.03445 \times 10^{-2}$   \\
		$3$ & $3.95 \times 10^{-4}$ & $1.01025 \times 10^{-2}$   \\
		$4$ & $1.98 \times 10^{-4}$ & $9.98198 \times 10^{-3}$   \\
		$5$ & $9.88 \times 10^{-5}$ & $9.92182 \times 10^{-3}$   \\
		$6$ & $4.94 \times 10^{-5}$ & $9.89128 \times 10^{-3}$   \\
		\hline
		\multicolumn{3}{c}{$h = 0.003125$} \\ \hline
		\multicolumn{1}{c}{$k$} & \multicolumn{1}{c}{$\Delta t$} & \multicolumn{1}{c}{$Error$}  \\
		\hline
		$0$ & $1.12 \times 10^{-3}$ & $1.94633 \times 10^{-3}$   \\
		$1$ & $5.59 \times 10^{-4}$ & $1.57935 \times 10^{-3}$   \\
		$2$ & $2.80 \times 10^{-4}$ & $1.39709 \times 10^{-3}$   \\
		$3$ & $1.40 \times 10^{-4}$ & $1.30627 \times 10^{-3}$   \\
		$4$ & $6.99 \times 10^{-5}$ & $1.26087 \times 10^{-3}$   \\
		$5$ & $3.49 \times 10^{-5}$ & $1.23823 \times 10^{-3}$   \\
		$6$ & $1.75 \times 10^{-5}$ & $1.22691 \times 10^{-3}$   \\
		\hline 
	\end{tabular}
\end{table}

\begin{table}[ht!]
	\centering
	\caption{Example~\ref{ex:1d} by {\rm (S2)}: reducing the time-step as  $\Delta t(k)=\frac{h}{2^{k}}$.}
	\label{tableX2}
	\begin{tabular}{rrr}
		\hline 
		\multicolumn{3}{c}{$h = 0.025$} \\ \hline
		\multicolumn{1}{c}{$k$} & \multicolumn{1}{c}{$\Delta t$} & \multicolumn{1}{c}{$Error$}  \\
		\hline
		$0$ & $2.50 \times 10^{-2}$ & $5.63 \times 10^{-3}$   \\
		$1$ & $1.25 \times 10^{-2}$ & $1.50 \times 10^{-3}$   \\
		$2$ & $6.25 \times 10^{-3}$ & $4.30 \times 10^{-4}$   \\
		$3$ & $3.13 \times 10^{-3}$ & $1.58 \times 10^{-4}$   \\
		$4$ & $1.56 \times 10^{-3}$ & $8.97 \times 10^{-5}$   \\
		$5$ & $7.81 \times 10^{-4}$ & $7.27 \times 10^{-5}$   \\
		$6$ & $3.91 \times 10^{-4}$ & $6.84 \times 10^{-5}$   \\
		\hline
		\multicolumn{3}{c}{$h = 0.003125$} \\ \hline
		\multicolumn{1}{c}{$k$} & \multicolumn{1}{c}{$\Delta t$} & \multicolumn{1}{c}{$Error$}  \\
		\hline
		$0$ & $3.13 \times 10^{-3}$ & $8.96 \times 10^{-5}$   \\
		$1$ & $1.56 \times 10^{-3}$ & $2.34 \times 10^{-5}$   \\
		$2$ & $7.81 \times 10^{-4}$ & $6.64 \times 10^{-6}$   \\
		$3$ & $3.91 \times 10^{-4}$ & $2.41 \times 10^{-6}$   \\
		$4$ & $1.95 \times 10^{-4}$ & $1.36 \times 10^{-6}$   \\
		$5$ & $9.77 \times 10^{-5}$ & $1.10 \times 10^{-6}$   \\
		$6$ & $4.88 \times 10^{-5}$ & $1.03 \times 10^{-6}$   \\
		\hline 
	\end{tabular}
\end{table}

\subsection{Examples for the two-dimensional case~$(d=2)$}
%
%
%
%
We set the next example in two-dimensional space.
\begin{Ex}[$d=2$]\label{ex:2d}
	In problem~\eqref{eqns:prob}, let $d=2$, $\Omega = (0,1)^d$ and $T=1$.
	We consider three functions for the velocity:
\begin{eqnarray*}
 & (i)~u(x,t) = (t,t)^{ \top }, \quad
 (ii)~u(x,t) = (x_1+t,x_2+t)^{ \top }, \\
 & (iii)~u(x,t) = (\sin(x_1+x_2+t), \sin(x_1+x_2+t))^{ \top },
\end{eqnarray*}
	which imply $\Gamma_\mathrm{in} = \{ (s,0)^{ \top }\in\pz\Omega;\ s\in [0, 1] \} \cup \{ (0,s)^{ \top }\in\pz\Omega;\ s\in [0, 1] \}$~$(t\in (0,T])$.
	The functions~$F$, $\zeta_\mathrm{in}$ and $\zeta^0$ are given so that the exact solution is
	\[
	 \zeta(x,t) = 
	 \begin{bmatrix}
	  \sin(x_1 + x_2 + t) + 2 & \sin(x_1 + x_2 + t) \\
	  \sin(x_1 + x_2 + t) & \sin(x_1 + x_2 + t) + 2
	 \end{bmatrix}.
	\]
\end{Ex}
\par
We solve Example~\ref{ex:2d} by (S1) with $\Delta t = c\sqrt{h}$ for $c = 1/20$ and (S2) with $\Delta t = c^\prime h$ for $c^\prime = 1/10$, where the mesh is constructed for $h_1=h_2=h=1/N$, i.e., $N_1=N_2=N$, with $N = 10, 20, 40$ and $80$.
Tables~\ref{table:ex2_s1_dt_sqrth} and \ref{table:ex2_s2_dt_h} show the values of error $E_{11}$ and their slopes in~$\Delta t$. Slope results for~$E_{12}$ and~$E_{22}$ adopting different velocity fields $(i)$, $(ii)$ and $(iii)$ are very similar to those obtained for $E_{11}$; thus they are omitted here in order to save space.
We can confirm that (S1) and (S2) are of second-order in $\Delta t$ in two-dimensional space for the three cases of velocity, $(i)$, $(ii)$ and $(iii)$.
These results are consistent with the theoretical results in Theorem~\ref{prop:sec_order_p_order}.
\begin{table}[!htbp]
	\centering
	\caption{{Example~\ref{ex:2d} by {\rm (S1)} with $\Delta t = c\sqrt{h}$~$(c=1/20)$:} Values of~$E_{11}$ and their slopes in $\Delta t$.}
	\begin{tabular}{rrrcrrcrr}
	\toprule
	& \multicolumn{2}{c}{$(i)$} && \multicolumn{2}{c}{$(ii)$} && \multicolumn{2}{c}{$(iii)$} \\ \cline{2-3}\cline{5-6}\cline{8-9}
	\multicolumn{1}{c}{$N$} & \multicolumn{1}{c}{$E_{11}$} & \multicolumn{1}{c}{Slope} && \multicolumn{1}{c}{$E_{11}$} & \multicolumn{1}{c}{Slope} && \multicolumn{1}{c}{$E_{11}$} & \multicolumn{1}{c}{Slope} \\
	\hline
	$10$ & $3.87 \times 10^{-2}$ & -- && $3.84 \times 10^{-2}$ & -- && $3.87 \times 10^{-2}$ & -- \\
	$20$ & $1.98 \times 10^{-2}$ & $1.94$ && $1.96 \times 10^{-2}$ & $1.94$ && $1.98 \times 10^{-2}$ & $1.94$  \\
	$40$ & $9.99 \times 10^{-3}$ & $1.97$ && $9.94 \times 10^{-3}$ & $1.97$ && $9.99 \times 10^{-3}$ & $1.97$ \\
	$80$ & $5.03 \times 10^{-3}$ & $1.98$ && $5.01 \times 10^{-3}$ & $1.98$ && $5.03 \times 10^{-3}$ & $1.98$    \\ 
	\bottomrule
	\end{tabular}
	\label{table:ex2_s1_dt_sqrth}
\end{table}
\begin{table}[!htbp]
	\centering
	\caption{{Example~\ref{ex:2d} by {\rm (S2)} with $\Delta t = c^\prime h$~$(c^\prime=1/10)$:} Values of~$E_{11}$ and their slopes in $\Delta t$.}
	\begin{tabular}{rrrcrrcrr}
	\toprule
	& \multicolumn{2}{c}{$(i)$} && \multicolumn{2}{c}{$(ii)$} && \multicolumn{2}{c}{$(iii)$} \\ \cline{2-3}\cline{5-6}\cline{8-9}
	\multicolumn{1}{c}{$N$} & \multicolumn{1}{c}{$E_{11}$} & \multicolumn{1}{c}{Slope} && \multicolumn{1}{c}{$E_{11}$} & \multicolumn{1}{c}{Slope} && \multicolumn{1}{c}{$E_{11}$} & \multicolumn{1}{c}{Slope} \\
	\hline
	$10$ & $2.07 \times 10^{-4}$ & -- && $2.18 \times 10^{-3}$ & -- && $9.79 \times 10^{-4}$ & -- \\
	$20$ & $5.10 \times 10^{-5}$ & $2.02$ && $5.35 \times 10^{-4}$ & $2.02$ && $2.53 \times 10^{-4}$ & $1.95$ \\
	$40$ & $1.27 \times 10^{-5}$ & $2.00$ && $1.32 \times 10^{-4}$ & $2.02$ && $6.39 \times 10^{-5}$ & $1.98$ \\
	$80$ & $3.17 \times 10^{-6}$ & $2.00$ && $3.27 \times 10^{-5}$ & $2.01$ && $1.61 \times 10^{-5}$ & $1.99$ \\
	\bottomrule
	\end{tabular}
	\label{table:ex2_s2_dt_h}
\end{table}
%
%
%
%
%
%
%
%
\subsection{The Oldroyd-B constitutive equation in two-dimensional space}
%
%
%
%
We apply our approximations of the upper-convected time derivative of second-order in time~\eqref{eq:sec_order_p_order} in Theorem~\ref{prop:sec_order_p_order} for solving a problem governed by the Oldroyd-B constitutive equation in two-dimensional space; find $\zeta: \Omega\times (0, T) \to \R^{d\times d}_\sym$ such that
	\begin{subequations}\label{eqns:prob_ob}
		\begin{align}
		&&&& \zeta + Wi \, \uctd{\zeta} & = 2 ( 1 - \beta ) D(u) + F && \mbox{in}\ \Omega\times (0, T), &&&&
		\label{eq:prob_ob1}\\
		&&&& \zeta & = \zeta_\mathrm{in} && \mbox{on}\ \Gamma_\mathrm{in}\times (0, T), &&&&
		\label{eq:prob_ob2}\\
		&&&& \zeta & = \zeta^0 && \mbox{in}\ \Omega,\ \mbox{at}\ t=0. &&&&
		\end{align}
	\end{subequations}
\par
The scheme to solve problem~\eqref{eqns:prob_ob} is to find $\{\zeta_h^n \in V_h;\ n=1, \ldots, N_T\}$ such that
\begin{subequations}\label{scheme_ob}
	\begin{align}
	\label{scheme_ob:general_step}
	\zeta_h^n(x) + Wi \, [\mathcal{A}_h^{n,(p)} \zeta_h] (x)
	& = 2(1-\beta)D(u^n)(x) + F^n(x), & x & \in \bar\Omega_h,\ n  \ge 1, \\
	\label{scheme_ob:initial_step}
	\zeta_h^0(x) & = \zeta^0(x), & x & \in \bar\Omega_h,
	\end{align}
\end{subequations}
where $\mathcal{A}_h^{n,(p)} \zeta_h:\bar\Omega_h\to\R^{d\times d}_\sym$ is the function defined already by~\eqref{def:An_h}.
When an upwind point is outside the domain, we employ a value of $\zeta_\mathrm{in}$ at closest lattice point to the upwind point similarly to the case of scheme~\eqref{scheme} as mentioned in Remark~\ref{rmk:upwind_cell_2d}-(iii).
In the following, scheme~\eqref{scheme_ob} with $p=1$ and $p=2$ for problem~\eqref{eqns:prob_ob} are called $({\rm S}1)^\prime$ and $({\rm S}2)^\prime$, respectively.
\par
We set two examples below:
\begin{Ex}[$d=2$]\label{ex:OldB2D}
	In problem~\eqref{eqns:prob_ob}, let $d=2$, $\Omega = (0,1)^d$, $T=1$ and $\beta = 1/9$.
	We consider six values of the {\rm Weissenberg} number $Wi$,
	\[
		Wi = 0.025, 1, 5, 10, 50, 100,
	\]
	and the following function for the velocity field:
\begin{align*}
 u(x,t) & = (\sin(x_1+x_2+t), \sin(x_1+x_2+t))^{ \top },
\end{align*}
	which implies $\Gamma_\mathrm{in} = \{ (s,0)^{ \top }\in\pz\Omega;\ s\in [0, 1] \} \cup \{ (0,s)^{ \top }\in\pz\Omega;\ s\in [0, 1] \}$.
	The functions~$F$, $\zeta_{\rm in}$ and $\zeta^0$ are given so that the exact solution is 
	\[
	 \zeta(x,t) = 
	 \begin{bmatrix}
	  \sin(x_1 + x_2 + t) + 2 & \sin(x_1 + x_2 + t) \\
	  \sin(x_1 + x_2 + t) & -\sin(x_1 + x_2 + t) + 2
	 \end{bmatrix}.
	\]
\end{Ex}
\begin{Ex}[$d=2$, {\cite{Venkatesan2017}}]\label{ex:Venkatesan2017}
	In problem~\eqref{eqns:prob_ob}, let $d=2$, $\Omega = (0,1)^d$, $T=0.5$, $\beta = 0.75$ and $Wi = 0.25$.
	We consider the following function for the velocity field:
\[
u(x, t) = (\exp(-0.1t)\sin(\pi x_1), -\pi\exp(-0.1t)x_2\cos(\pi x_1))^{ \top },
\]
	which implies $\Gamma_{\rm in} = \{ (s,0)^{ \top }\in\pz\Omega;\ s\in [0, 1] \} \cup \{ (0,s)^{ \top }\in\pz\Omega;\ s\in [0, 1] \}$.
	The functions~$F$, $\zeta_{\rm in}$ and $\zeta^0$ are given so that the exact solution is
\[
\zeta(x,t) = 
\begin{bmatrix}
	\exp(-0.1t)\sin(\pi x_1) & -\pi\exp(-0.1t)x_2\cos(\pi x_1) \\ 
	-\pi\exp(-0.1t)x_2\cos(\pi x_1) & \exp(-0.1t)\sin(\pi x_1)\cos(\pi x_2) \end{bmatrix}.
\]
\end{Ex}
\par
We solve Example~\ref{ex:OldB2D} by $({\rm S}1)^\prime$ with $\Delta t = c\sqrt{h}$ for $c = 1/50$ and $({\rm S}2)^\prime$ with $\Delta t = c^\prime h$ for $c^\prime = 1/5$, where the mesh is constructed for $h_1=h_2=h=1/N$, i.e., $N_1=N_2=N$, with $N = 10, 20, 40$ and $80$.
In order to further investigate the errors and the orders of convergence of the schemes for solving problem~\eqref{eqns:prob_ob}, we give the results for the three different components $\zeta_{11}$, $\zeta_{12}$ and $\zeta_{22}$.
Tables~\ref{table:ErrorOB2D1stOrder} and~\ref{table:ErrorOB2D2stOrder} show the results by $({\rm S}1)^\prime$ and $({\rm S}2)^\prime$, respectively, for~$Wi = 0.025$.
From a quantitative point of view, the results are consistent with the theoretical results in Theorem~\ref{prop:sec_order_p_order}.
\begin{table}[htbp!]
		\centering
		\caption{Example~\ref{ex:OldB2D} by $({\rm S}1)^\prime$ with $\Delta t = c\sqrt{h}$~$(c=1/50)$: Values of each tensor entry ~$E_{11},E_{12},E_{22}$ and their slopes in $\Delta t$ for $Wi = 0.025$.}
		\begin{tabular}{rrrcrrcrr}
			\toprule
			\multicolumn{1}{c}{$N$} & \multicolumn{1}{c}{$E_{11}$} & \multicolumn{1}{c}{Slope} && \multicolumn{1}{c}{$E_{12}$} & \multicolumn{1}{c}{Slope} && \multicolumn{1}{c}{$E_{22}$} & \multicolumn{1}{c}{Slope} \\
		\hline
			$10$ & $2.03 \times 10^{-3}$ & -- && $2.03 \times 10^{-3}$ & -- && $2.03 \times 10^{-3}$ & -- \\
			$20$ & $1.02 \times 10^{-3}$ & $1.99 $ && $1.02 \times 10^{-3}$ & $1.99 $ && $ 1.02 \times 10^{-3}$ & $1.99 $ \\
			$40$ & $5.11 \times 10^{-4}$ & $ 1.99$ && $5.11 \times 10^{-4}$ & $1.99 $ && $ 5.11 \times 10^{-4}$ & $ 1.99$ \\
			$80$ & $2.56 \times 10^{-4}$ & $1.99 $ && $2.56 \times 10^{-4}$ & $1.99 $ && $2.56 \times 10^{-4}$ & $ 1.99$ \\
			\bottomrule
	\end{tabular}
	\label{table:ErrorOB2D1stOrder}
\end{table}
\begin{table}[htbp!]
	\centering
	\caption{Example~\ref{ex:OldB2D} by $({\rm S}2)^\prime$ with $\Delta t = c^\prime h$~$(c^\prime=1/5)$: Values of each tensor entry~$E_{11},E_{12},E_{22}$ and their slopes in $\Delta t$ for $Wi = 0.025$.}
	\begin{tabular}{rrrcrrcrr}
		\toprule
		\multicolumn{1}{c}{$N$} & \multicolumn{1}{c}{$E_{11}$} & \multicolumn{1}{c}{Slope} && \multicolumn{1}{c}{$E_{12}$} & \multicolumn{1}{c}{Slope} && \multicolumn{1}{c}{$E_{22}$} & \multicolumn{1}{c}{Slope} \\
		\hline
		$10$ & $7.62 \times 10^{-5}$ & -- && $7.24 \times 10^{-5}$ & -- && $7.62 \times 10^{-5}$ & -- \\
		$20$ & $1.89 \times 10^{-6}$ & $2.02 $ && $1.80 \times 10^{-6}$ & $2.01 $ && $1.89 \times 10^{-5}$ & $2.02 $ \\
		$40$ & $4.75 \times 10^{-6}$ & $1.99 $ && $4.57 \times 10^{-6}$ & $ 1.98$ && $4.75 \times 10^{-6}$ & $ 1.99$ \\
		$80$ & $1.21 \times 10^{-6}$ & $ 1.97$ && $1.17 \times 10^{-6}$ & $ 1.96$ && $ 1.21 \times 10^{-6}$ & $1.97 $ \\
		\bottomrule
	\end{tabular}
	\label{table:ErrorOB2D2stOrder}
\end{table}
\par

A computational challenge in viscoelastic fluid flows is the application of high values of the Weissenberg number, i.e., $Wi>1$. In fact, the \textit{infamous} High Weissenberg Number Problem \cite{Fattal:2004,Hulsen,Martins2015} depends of some particular factors on viscoelastic flows, as for instance, domain geometry, boundary conditons, fluid type, mesh size, etc. In summary, this instability is related to the unbounded values of the stress tensor during the transient solution resulting in the fail of the numerical methods. It is important to highlight that some classical methods, i.e. without stabilization techniques, have failed for $Wi= O(1)$ exhibiting numerical oscillations of the solution. Rougly speaking, the High Weissenberg Number Problem can be interpreted according to a critical value of the Weissenberg number, $Wi_{crit}$, which the numerical solution is boundly maintained during the simulation of the classical constitutive formulations. For example, considering the traditional Oldroyd-B model, Fattal and Kupferman \cite{Fattal:2004} described  $Wi_{crit} \approx 0.5$ for the cavity flow while Oliveira and Miranda \cite{OliveiraMiranda} pointed out $Wi_{crit} \approx 1$ for unsteady viscoelastic flow past bounded cylinders. Moreover, Walters and Webster \cite{Walters2003} presented results for the $4:1$ contraction problem with the critical Weissenberg number near to 3. Therefore, there is an effort of the researchers to circumvent the High Weissenberg Number Problem developing new formulations that can be stable in simulations with $Wi > Wi_{crit}$.

It is important to highlight that the schemes presented in this current work can deal with high values of $Wi$ without the need to employ stabilization strategies. To test the accuracy of $({\rm S}2)^\prime$, we vary the values of Weissenberg number as $Wi=1,5,10,50,100$ in Example~\ref{ex:OldB2D} and the results are presented in Table \ref{table:ErrorOB2D2nd_Wi_varied1}. The main focus for varying the Weissenberg number is to verify the ability of $({\rm S}2)^\prime$ for dealing with the Oldroyd-B constitutive equation defined on the context of high elasticity. From the results presented in Table \ref{table:ErrorOB2D2nd_Wi_varied1}, we can notice that the numerical order of convergence of $({\rm S}2)^\prime$ is of second-order in both time and space, and that the effect of varying the Weissenberg number is not significant for this example.

	\begin{table}[ht!]
	\centering
	\caption{Example~\ref{ex:OldB2D} by $({\rm S}2)^\prime$ with $\Delta t = c^\prime h$~$(c^\prime=1/5)$ and different values of $Wi$ number.}
	\label{table:ErrorOB2D2nd_Wi_varied1}
	\begin{tabular}{rrrcrrcrr}
		\hline 
		\multicolumn{9}{c}{$Wi=1.0$} \\ \hline
		\multicolumn{1}{c}{$N$} & \multicolumn{1}{c}{$E_{11}$} & \multicolumn{1}{c}{Slope} && \multicolumn{1}{c}{$E_{12}$} & \multicolumn{1}{c}{Slope} && \multicolumn{1}{c}{$E_{22}$} & \multicolumn{1}{c}{Slope} \\ \hline
		$10$ & $1.55 \times 10^{-3}$ & -- && $1.06 \times 10^{-3}$ & -- && $5.54 \times 10^{-4}$ & -- \\
		$20$ & $4.23 \times 10^{-4}$ & $1.88$ && $2.93 \times 10^{-4}$ & $1.85$ && $1.48 \times 10^{-4}$ & $1.91$ \\
		$40$ & $1.09 \times 10^{-4}$ & $1.95$ && $7.65 \times 10^{-5}$ & $1.94$ && $3.79 \times 10^{-5}$ & $1.96$ \\
		$80$ & $2.77 \times 10^{-5}$ & $1.98$ && $1.95 \times 10^{-5}$ & $1.98$ && $9.58 \times 10^{-6}$ & $1.99$ \\
		\hline
		\multicolumn{9}{c}{$Wi=5$} \\ \hline
		\multicolumn{1}{c}{$N$} & \multicolumn{1}{c}{$E_{11}$} & \multicolumn{1}{c}{Slope} && \multicolumn{1}{c}{$E_{12}$} & \multicolumn{1}{c}{Slope} && \multicolumn{1}{c}{$E_{22}$} & \multicolumn{1}{c}{Slope} \\ \hline
		$10$ & $1.97 \times 10^{-3}$ & -- && $1.37 \times 10^{-3}$ & -- && $7.13 \times 10^{-4}$ & -- \\
		$20$ & $5.36 \times 10^{-4}$ & $1.87$ && $3.80 \times 10^{-4}$ & $1.85$ && $1.97 \times 10^{-4}$ & $1.86$ \\
		$40$ & $1.39 \times 10^{-4}$ & $1.95$ && $9.90 \times 10^{-5}$ & $1.94$ && $5.14 \times 10^{-5}$ & $1.94$ \\
		$80$ & $3.51 \times 10^{-5}$ & $1.98$ && $2.52 \times 10^{-5}$ & $1.98$ && $1.31 \times 10^{-5}$ & $1.97$ \\
		\hline
		\multicolumn{9}{c}{$Wi=10$} \\ \hline
		\multicolumn{1}{c}{$N$} & \multicolumn{1}{c}{$E_{11}$} & \multicolumn{1}{c}{Slope} && \multicolumn{1}{c}{$E_{12}$} & \multicolumn{1}{c}{Slope} && \multicolumn{1}{c}{$E_{22}$} & \multicolumn{1}{c}{Slope} \\ \hline
		$10$ & $2.03 \times 10^{-3}$ & -- && $1.42 \times 10^{-3}$ & -- && $7.38 \times 10^{-4}$ & -- \\
		$20$ & $5.54 \times 10^{-4}$ & $1.87$ && $3.93 \times 10^{-4}$ & $1.85$ && $2.04 \times 10^{-4}$ & $1.85$ \\
		$40$ & $1.43 \times 10^{-4}$ & $1.94$ && $1.03 \times 10^{-4}$ & $1.94$ && $5.35 \times 10^{-5}$ & $1.93$ \\
		$80$ & $3.63 \times 10^{-5}$ & $1.98$ && $2.61 \times 10^{-5}$ & $1.98$ && $1.36 \times 10^{-5}$ & $1.97$ \\
		\hline
		\multicolumn{9}{c}{$Wi=50$} \\ \hline
		\multicolumn{1}{c}{$N$} & \multicolumn{1}{c}{$E_{11}$} & \multicolumn{1}{c}{Slope} && \multicolumn{1}{c}{$E_{12}$} & \multicolumn{1}{c}{Slope} && \multicolumn{1}{c}{$E_{22}$} & \multicolumn{1}{c}{Slope} \\ \hline
		$10$ & $2.08 \times 10^{-3}$ & -- && $1.46 \times 10^{-3}$ & -- && $7.59 \times 10^{-4}$ & -- \\
		$20$ & $5.69 \times 10^{-4}$ & $1.87$ && $4.05 \times 10^{-4}$ & $1.85$ && $2.11 \times 10^{-4}$ & $1.85$ \\
		$40$ & $1.47 \times 10^{-4}$ & $1.95$ && $1.06 \times 10^{-4}$ & $1.94$ && $5.53 \times 10^{-5}$ & $1.93$ \\
		$80$ & $3.72 \times 10^{-5}$ & $1.98$ && $2.68 \times 10^{-5}$ & $1.98$ && $1.41 \times 10^{-5}$ & $1.97$ \\
		\hline
		\multicolumn{9}{c}{$Wi=100$} \\ \hline
		\multicolumn{1}{c}{$N$} & \multicolumn{1}{c}{$E_{11}$} & \multicolumn{1}{c}{Slope} && \multicolumn{1}{c}{$E_{12}$} & \multicolumn{1}{c}{Slope} && \multicolumn{1}{c}{$E_{22}$} & \multicolumn{1}{c}{Slope} \\ \hline
		$10$ & $2.09 \times 10^{-3}$ & -- && $1.46 \times 10^{-3}$ & -- && $7.62 \times 10^{-4}$ & -- \\
		$20$ & $5.71 \times 10^{-4}$ & $1.87$ && $4.06 \times 10^{-4}$ & $1.85$ && $2.12 \times 10^{-4}$ & $1.85$ \\
		$40$ & $1.48 \times 10^{-4}$ & $1.95$ && $1.06 \times 10^{-4}$ & $1.94$ && $5.55 \times 10^{-5}$ & $1.93$ \\
		$80$ & $3.74 \times 10^{-5}$ & $1.98$ && $2.69 \times 10^{-5}$ & $1.98$ && $1.42 \times 10^{-5}$ & $1.97$ \\
		\hline 
	\end{tabular}
\end{table}
\par
Finally, example~\ref{ex:Venkatesan2017} employs the manufactured solution used by Venkatesan and Ganesan~\cite{Venkatesan2017}.
Notice that in this study we are investigating the numerical behavior of the schemes for non-homogeneous boundary conditions in parts of the domain.
Table~\ref{table:ErrorOB2_Venkatesan2017} describes the results for Example~\ref{ex:Venkatesan2017} by $({\rm S}2)^\prime$ with $\Delta t = c^\prime h$ for $c^\prime = 1/10$, where the mesh is constructed for $h_1=h_2=h=1/N$, i.e., $N_1=N_2=N$, with $N = 10, 20, 40$ and $80$.
From this table we can see that these results are consistent with our truncation error analysis in Theorem~\ref{prop:sec_order_p_order}.
\begin{table}[htbp!]
	\centering
	\caption{Example~\ref{ex:Venkatesan2017} by $({\rm S}2)^\prime$ with $\Delta t = c^\prime h$~$(c^\prime=1/10)$: Values of each tensor entry~$E_{11},E_{12},E_{22}$ and their slopes in $\Delta t$ for $Wi = 0.25$ and $\beta = 0.75$.}
	\begin{tabular}{rrrcrrcrr}
		\toprule
		\multicolumn{1}{c}{$N$} & \multicolumn{1}{c}{$E_{11}$} & \multicolumn{1}{c}{Slope} && \multicolumn{1}{c}{$E_{12}$} & \multicolumn{1}{c}{Slope} && \multicolumn{1}{c}{$E_{22}$} & \multicolumn{1}{c}{Slope} \\
		\hline
		$10$ & $4.10 \times 10^{-3}$ & -- && $7.64 \times 10^{-2}$ & -- && $1.98 \times 10^{-2}$ & -- \\
		$20$ & $1.02 \times 10^{-3}$ & $2.01$ && $2.11 \times 10^{-3}$ & $1.86$ && $5.19 \times 10^{-3}$ & $1.93$ \\
		$40$ & $2.82 \times 10^{-4}$ & $1.86$ && $5.83 \times 10^{-4}$ & $1.85$ && $1.32 \times 10^{-3}$ & $1.97$ \\
		$80$ & $7.47 \times 10^{-5}$ & $1.91$ && $1.54 \times 10^{-4}$ & $1.92$ && $3.30 \times 10^{-4}$ & $2.00$ \\
		\bottomrule
	\end{tabular}
	\label{table:ErrorOB2_Venkatesan2017}
\end{table}
%
%
%
%
\section{Conclusions}
\label{sec:conclusions}
%
%
%
%
The application of the generalized Lie derivative~(GLD) for constructing schemes to deal with the upper-convected time derivative is an alternative form in the numerical solution of constitutive equations.
In spite of the success of this strategy firstly proposed by Lee and Xu~\cite{Lee2006}, to the best knowledge of the authors, the methodology was only applied in the context of finite elements.
In this work, we have combined a Lagrangian framework with GLD to develop new second-order finite difference approximations for the upper-convected time derivative.
Particularly, the schemes are constructed based on bilinear and biquadratic interpolation operators for solving a simple model in one- and two-dimensional spaces.
The schemes are explicit and no CFL condition is required as the Lagrangian framework is employed.
Truncation errors of~$O(\Delta t^2 + h^p)$ $(p=1,2)$ for the finite difference approximations of the the upper-convected time derivative have been proved.
A numerical integration of composite functions may cause an instability in the case of Lagrangian finite element method, our schemes, however, do not have such instability since there is no numerical integration thanks to the finite difference method.
According to our numerical results for simplified model equations, the new finite difference schemes can reach second-order of accuracy in time and {space~$(p=2)$} corroborating with the theoretical analysis. Moreover, the proposed strategy has been also applied to solve a two-dimensional Oldroyd-B constitutive equation subjected to a prescribed velocity field. The results have been very satisfactory since the increasing of the Weissenberg number did not influence the good properties of accuracy and stability of the finite difference approximations. As a future work, we intend to extend our schemes for solving viscoelastic fluid flows governed by different constitutive equations at high Weissenberg numbers.
\appendix
%
\section*{Appendix}
\renewcommand{\thesection}{A}
\setcounter{Lem}{0}
\renewcommand{\theLem}{\thesection.\arabic{Lemma}}
\setcounter{Rmk}{0}
\renewcommand{\theRmk}{\thesection.\arabic{Remark}}
\setcounter{figure}{0}
\renewcommand{\thefigure}{\thesection.\arabic{figure}}
\renewcommand{\thealgorithm}{\thesection.\arabic{algorithm}}
\setcounter{equation}{0}
\makeatletter
  \renewcommand{\theequation}{%
  \thesection.\arabic{equation}}
  \@addtoreset{equation}{section}
\makeatother
%
%
%
\subsection{Proofs of properties in~\eqref{eqns:L}}\label{A.subsec:L}
%
%
Firstly, we prove~\eqref{eqns:L1}.
The second equality of~\eqref{eqns:L1} is obtained immediately from the definition of~$L$ in~\eqref{def:L} as
\begin{align*}
L_{ij} (x,t;t_1,t_1)
& = \Bigl[ \prz{}{z_j} X_i(z,t_1; t_1) \Bigr]_{{\dps |}z=X(x,t; t_1)}
= \Bigl[ \prz{}{z_j} z_i \Bigr]_{{\dps |}z=X(x,t; t_1)} = \bigl[ \delta_{ij} \bigr]_{{\dps |}z=X(x,t; t_1)}
= \delta_{ij},
\end{align*}
where $\delta_{ij}$~$(i,j=1,\dots,d)$ is Kronecker's delta function.
For the first equality of~\eqref{eqns:L1}, we prove
\begin{equation}
I = L(x,t;t_1,t_2) L(x,t;t_2,t_1).
\label{eq:I_LL}
\end{equation}
Let $x\in\bar\Omega$ and $t_1, t_2 \in [0,T]$ be fixed arbitrarily.
For any $y\in\bar\Omega$, it holds that
\[
y = X \bigl( X(y,t_2; t_1),t_1; t_2 \bigr),
\]
which is equivalent to
\begin{equation}
y_i = X_i(X(y,t_2; t_1),t_1; t_2), \quad i=1, \ldots, d.
\label{eq:phi_identity}
\end{equation}
The differentiation of both sides of~\eqref{eq:phi_identity} with respect to~$y_j~(j=1,\ldots,d)$ implies that
\begin{align}
\delta_{ij}
& = \prz{}{y_j} \Bigl( X_i(X(y,t_2; t_1),t_1; t_2) \Bigr) 
= \sum_{k=1}^d \Bigl[ \prz{}{z_k} X_i(z,t_1; t_2) \Bigr]_{\dps |z = X(y,t_2; t_1)} \Bigl[ \prz{}{z_j} X_k(z,t_2; t_1) \Bigr]_{\dps | z = y}. \label{eq:proof_L1}
\end{align}
Substituting $X(x,t; t_2)$ into $y$ in~\eqref{eq:proof_L1} and using $X( X(x,t; t_2), t_2; t_1) = X(x,t; t_1)$, we get
\begin{align*}
\delta_{ij} & = \sum_{k=1}^d \Bigl[ \prz{}{z_k} X_i(z,t_1; t_2) \Bigr]_{\dps |z = X(x,t; t_1)} \Bigl[ \prz{}{z_j} X_k(z,t_2; t_1) \Bigr]_{\dps | z = X(x,t; t_2)} 
= \sum_{k=1}^d L_{ik}(x,t;t_1,t_2) L_{kj}(x,t;t_2,t_1),
\end{align*}
which implies~\eqref{eq:I_LL}.
Thus, the first equality of~\eqref{eqns:L1} holds true.
\par
Secondly, we prove~\eqref{eqns:L2}.
From the definition of $L$ in~\eqref{def:L}, we have
\begin{align*}
& \prz{}{s} L_{ij}(x,t;t_1,s) 
= \Bigl[ \prz{}{s} \prz{}{z_j} X_i(z,t_1; s) \Bigr]_{\dps |z = X(x,t; t_1)} \\
& = \Bigl[ \prz{}{z_j} \prz{}{s} X_i(z,t_1; s) \Bigr]_{\dps |z = X(x,t; t_1)} 
= \Bigl[ \prz{}{z_j} u_i \bigl( X(z,t_1; s), s \bigr) \Bigr]_{\dps |z = X(x,t; t_1)} \\
& = \biggl[ \sum_{k=1}^d \prz{u_{i}}{x_k} \bigl( X(z,t_1; s), s \bigr) 
\prz{}{z_j} X_k(z,t_1; s) \biggr]_{\dps | z = X(x,t; t_1)} 
= \sum_{k=1}^d \prz{u_{i}}{x_k} \bigl( X(x,t; s), s \bigr) \Bigl[ \prz{}{z_j} X_k(z,t_1; s) \Bigr]_{\dps | z = X(x,t; t_1)} \\
& = \sum_{k=1}^d [\nabla u]_{ik} \bigl( X(x,t; s), s \bigr) L_{kj}(x,t;t_1,s),
\end{align*}
which implies~\eqref{eqns:L2}.
\par
Finally, we prove~\eqref{eqns:L3}.
Property~\eqref{eqns:L1} gives an identity,
\begin{displaymath}
I = L(x,t;t_1,s) L(x,t;s,t_1).
\end{displaymath}
Considering the derivative of the identity above with respect to $s$, we have
\begin{align*}
0
& = \prz{}{s} \Bigl[ L(x,t;t_1,s) L(x,t;s,t_1) \Bigr] 
= \Bigl[ \prz{}{s} L(x,t;t_1,s) \Bigr] L(x,t;s,t_1) + L(x,t;t_1,s) \Bigl[ \prz{}{s} L(x,t;s,t_1) \Bigr] \\
& = (\nabla u) \bigl( X(x,t; s), s \bigr) L(x,t;t_1,s) L(x,t;s,t_1) + L(x,t;t_1,s) \Bigl[ \prz{}{s} L(x,t;s,t_1) \Bigr] \qquad \mbox{(by~\eqref{eqns:L2})} \\
& = (\nabla u) \bigl( X(x,t; s), s \bigr) + L(x,t;t_1,s) \Bigl[ \prz{}{s} L(x,t;s,t_1) \Bigr] \qquad \mbox{(by~\eqref{eqns:L1})},
\end{align*}
which completes the proof of~\eqref{eqns:L3} as
\begin{align*}
\prz{}{s} L(x,t;s,t_1)
& = - L (x,t;t_1,s)^{-1} (\nabla u) \bigl( X(x,t; s), s \bigr) 
= - L (x,t;s,t_1) (\nabla u) \bigl( X(x,t; s), s \bigr) \qquad \mbox{(by~\eqref{eqns:L1}).}
\end{align*}
%
%
%
\subsection{Proof of~\eqref{eq:UCM_GLD}}\label{A.subsec:UCM_GLD}
%
For the sake of simplicity, we employ simple notations, $L(\cdot,\cdot) = L(x,t;\,\cdot\,, \cdot\,)$ and $X = X(x,t; \,\cdot\,)$, as there is no confusion.
From the definition of the generalized Lie derivative in~\eqref{def:GLD} and the properties of $L$ in~\eqref{eqns:L}, we have
\begin{align*}
& (\mcL_u\zeta) (x,t) 
= (\mcL_u\zeta)\bigl( X(s), s \bigr)_{{\dps |}s=t} 
= L(t,s) \prz{}{s}\Bigl[ L(s,t) \zeta \bigl( X(s), s \bigr) L(s,t)^{ \top } \Bigr] L(t,s)^{ \top } {}_{{\dps |}s=t} \notag\\
& = L(t,s) \biggl[
\Bigl( \prz{}{s} L(s,t) \Bigr) \zeta \bigl( X(s), s \bigr) L(s,t)^{ \top } 
+ L(s,t) \Bigl( \prz{}{s} \zeta \bigl( X(s), s \bigr) \Bigr) L(s,t)^{ \top } 
+ L(s,t) \zeta \bigl( X(s), s \bigr) \Bigl( \prz{}{s} L(s,t)^{ \top } \Bigr)
\biggr] L(t,s)^{ \top } {}_{{\dps |}s=t} \notag\\
& = L(t,s) \biggl[
\Bigl( - L(s,t) (\nabla u) \bigl( X(s), s \bigr) \Bigr) \zeta \bigl( X(s), s \bigr) L(s,t)^{ \top }
+ L(s,t) \fz{D\zeta}{Dt} \bigl( X(s), s \bigr) L(s,t)^{ \top } \notag\\
& \qquad + L(s,t) \zeta \bigl( X(s), s \bigr) \Bigl( - L(s,t) (\nabla u) \bigl( X(s), s \bigr) \Bigr)^{ \top }
\biggr] L(t,s)^{ \top } {}_{{\dps |}s=t} \notag\\
& = L(t,s) \biggl[
- L(s,t) (\nabla u) \bigl( X(s), s \bigr) \zeta \bigl( X(s), s \bigr) L(s,t)^{ \top } 
+ L(s,t) \fz{D\zeta}{Dt} \bigl( X(s), s \bigr) L(s,t)^{ \top } \notag\\
& \qquad - L(s,t) \zeta \bigl( X(s), s \bigr) (\nabla u)^{ \top } \bigl( X(s), s \bigr) L(s,t)^{ \top } 
\biggr] L(t,s)^{ \top } {}_{{\dps |}s=t} \notag\\
& = \biggl[
- (\nabla u) \bigl( X(s), s \bigr) \zeta \bigl( X(s), s \bigr)
+ \fz{D\zeta}{Dt} \bigl( X(s), s \bigr) 
- \zeta \bigl( X(s), s \bigr) (\nabla u)^{ \top } \bigl( X(s), s \bigr)
\biggr]_{{\dps |}s=t} \notag\\
& = - (\nabla u) (x,t) \zeta (x,t)
+ \fz{D\zeta}{Dt} (x,t) 
- \zeta (x,t) (\nabla u)^{ \top } (x,t), \notag
\end{align*}
which completes the proof of~\eqref{eq:UCM_GLD}.
%
%
%
%
\subsection{Pseudo-codes for the proposed scheme}\label{A.subsec:PseudoCode}
%
The Algorithm~\ref{alg:Interpolation} contains the steps of the interpolation process for the evaluated function on the characteristic curve at a previous time.
\par
The main algorithm (see Algorithm~\ref{alg:General}) has all declarations and computations used to update the numerical solution on time.

\begin{algorithm}[htbp]
	\caption{Interpolation algorithm}
	\label{alg:Interpolation}
	\begin{algorithmic}[1]
		\Require ${y}^n_{i,j}$,$h_{d}(d=1,2)$, $(i,j)$, $\Lambda_\Omega$, $p$, $x_{i,j} $, $\bar{\Omega}_h$ and $\zeta_h$.
		\State Calculate the index on the discretized mesh $(i_0,j_0) \in \Lambda_\Omega$, i.e., $index(y^n_{i,j})=(i_0,j_0)$.
		\State Calculate $\eta_{i_0}^{(p)}$ and $\eta_{j_0}^{(p)}$ 
		\If{ $ p = 1$ }
		\begin{align*}
		\dps
		\eta_{i_0}^{(1)}(x; h_1) & \defeq
		\left\{
		\begin{aligned}
		& \frac{x - x_{{i_0}-1}}{h_1} && \bigl( x\in [x_{i_0-1}, x_{i_0} ) \bigr),\\
		& \frac{x_{i_0+1} - x}{h_1} && \bigl( x\in [x_{i_0}, x_{i_0 +1}] \bigr),\\
		& \ \ \ \ \   0 && \bigl( {\rm otherwise} \bigr),
		\end{aligned}
		\right. 
		\end{align*}
		and the same to calculate the function $\eta_{j_0}^{(1)}$ using the index $j_0$ and space-step $h_1$.
		\Else  
		\begin{align*}
		\dps
		\eta_{i_0}^{(2)}(x; h_1) & \defeq
		\left\{
		\begin{aligned}
		& \frac{x - x_{{i_0}-1}}{h_1}\cdot\frac{x - x_{i_0-2}}{2h_1} && \bigl( x\in [x_{i_0-2}, x_{i_0} ) \bigr),\\
		& \frac{x_{i_0+1} - x}{h_1}\cdot\frac{x - x_{i_0+2} }{2h_1} && \bigl( x\in [x_{i_0}, x_{i_0 +2}] \bigr),\\
		& \ \ \ \ \   0 && \bigl( {\rm otherwise} \bigr),
		\end{aligned}
		\right. 
		\end{align*}
		and the again to compute the function $\eta_{j_0}^{(2)}$ using the index $j_0$ and space-step $h_2$.
		\EndIf
		\State Define the basis function $\varphi_{i_0,j_0}^{(p)}$ as
		\begin{displaymath}
		\dps
		\varphi_{i_0,j_0}^{(p)}(y^n_{i,j}) \defeq 
		\eta_{i_0}^{(p)}(y^n_i;h_1) \eta_{j_0}^{(p)}(y^n_j;h_2).
		\end{displaymath}
		
		\State Compute the interpolation of the given function $\zeta_h$ at ${y}^n_{i,j}$ by 
		\[
		{Z}_{i,j}^{1,(p)}\defeq\bigl( \Pi_h^{(p)} \zeta_h \bigr) ({y}_{i,j}^n) = \sum_{x_{i,j}\in \bar\Omega_h} \zeta_h(x_{i,j})\varphi_{i_0,j_0}^{(p)}({y}_{i,j}^n).
		\]
		\Return ${Z}_{i,j}^{1,(p)} $. 
		
	\end{algorithmic}
\end{algorithm}%

\begin{algorithm}[htbp]
	\caption{Main algorithm}
	\label{alg:General}
	\begin{algorithmic}[1]
		\Require The domain $\Omega$ with $a_d~(d=1,2)$, division numbers $N_d~(d=1,2)$, interpolation order $p$, final time $T$, time step $\Delta t$ and the exact solution $\zeta_{\rm exact}^n(x_{i,j})$.
		\State Calculate $h_d=a_d/N_d~(d=1,2)$ and $\bar{\Omega}_h$, where $N_d~(d=1,2)$ are even numbers and $M_d = N_d/2~(d=1,2)$ for $p=2$, the indexes domain $\Lambda_\Omega = \{(i,j);\ i = 0,\ldots, N_1,\ j=0,\ldots, N_2 \}$ and the number of time steps $N_T$.
		\State Initialize the value $\zeta^n_{\rm in}(x_{i,j}) = \zeta^n_{\rm exact}(x_{i,j})$ for $n =0$ with 
		$\ (i,j)\in\Lambda_\Omega$ and $x_{i,j} \in \Gamma_\mathrm{in} = \{ (s,0)^\top\in\pz\Omega;\ s\in [0, 1] \} \cup \{ (0,s)^\top\in\pz\Omega;\ s\in [0, 1] \}$~$(t\in (0,T])$.
		\State Define the functions $u_{i,j}^n$, $\nabla u_{i,j}^n$ and $F^n_{i,j}$.
		\State Set $n=1$.
		\For{ $(i,j) \in \Lambda_\Omega$}
		\State Calculate the interpolation point $y^1_{i,j} \defeq X_1^{n} (x_{i,j}) = x_{i,j} - u^1_{i,j}\Delta t$.
		\If{ ($y^1_{i,j}\notin \bar{\Omega}$) }
		\[
		\zeta_{i,j}^1 = \zeta^1_{\rm in}(x_{i,j}).
		\]
		\Else  
		\State Compute $\nabla u_{i,j}^n$ and $F^n_{i,j}$ using the velocity field $u_{i,j}^n$. 
		\State Use Algorithm~\ref{alg:Interpolation} for ${y}^n_{i,j}$ to get ${Z}_{i,j}^{1,(p)}=\bigl( \Pi_h^{(p)} \zeta_h \bigr) ({y}_{i,j}^1)$.
		\State Update $\zeta^1_{i,j}$ by an approximation of first order in time
		\[
		\zeta_{i,j}^1 =
		\bigl[I + \Delta t (\nabla u^1_{i,j}) \bigr] Z_{i,j}^{1,(p)} \bigl[I + \Delta t (\nabla u^1_{i,j}) \bigr]^\top + \Delta t F^1_{i,j}.
		\]
		\EndIf
		\EndFor
		\While{ $n \le N_T$}
		\For{ $(i,j) \in \Lambda_\Omega$}
		\State Calculate the interpolation points $y^n_{i,j} = X_1^{n} (x_{i,j}) = x_{i,j} - u^n_{i,j}\Delta t$ and $\tilde{y}^n_{i,j} = \tilde{X}_1^{n} (x_{i,j}) = x_{i,j} - 2u^n_{i,j}\Delta t$.
		\If{ ($y^n_{i,j}\notin \bar{\Omega}$ or $\tilde{y}^n_{i,j}\notin \bar{\Omega}$) }
		\[
		\zeta_{i,j}^n = \zeta^n_{\rm in}(x_{i,j})
		\]
		\Else  
		\State Compute $\nabla u_{i,j}^n$ and $F^n_{i,j}$ using the velocity field $u_{i,j}^n$ 
		\State Use Algorithm~\ref{alg:Interpolation} for ${y}^n_{i,j}$ to get ${Z}_{i,j}^{n,(p)}=\bigl( \Pi_h^{(p)} \zeta_h \bigr) ({y}_{i,j}^n) $ and solve again for $\tilde{y}^n_{i,j}$ to get $\tilde{Z}_{i,j}^{n,(p)}=\bigl( \Pi_h^{(p)} \zeta_h \bigr) (\tilde{y}_{i,j}^n) $.
		\State Update $\zeta^n_{i.j}$ by an approximation of second order in time 
		\begin{align*}
		\zeta_{i,j}^n & =
		\frac{4}{3} \bigl[I + \Delta t (\nabla u^n_{i,j}) \bigr] Z_{i,j}^{n,(p)} \bigl[I + \Delta t (\nabla u^n_{i,j}) \bigr]^\top \\
		& \quad - \frac{1}{3} \bigl[I + 2\Delta t (\nabla u^n_{i,j}) \bigr] \tilde{Z}_{i,j}^{n,(p)} \bigl[I + 2\Delta t (\nabla u^n_{i,j}) \bigr]^\top + \frac{2\Delta t}{3} F^n_{i,j}.
		\end{align*}
		\EndIf
		\EndFor
		\State $n \leftarrow n+1$. 
		\EndWhile 
		
	\end{algorithmic}
\end{algorithm}

\clearpage

\section*{Acknowledgments}
{D.O.M.} would like to acknowledge the support of the grant 2013/07375-0 - Center of Mathematical Sciences Applied to industry (Cepid-CeMEAI), grants 2019/08742-2 and 2017/11428-2, S{\~a}o Paulo Research Foundation (FAPESP). {H.N.} is supported by JSPS KAKENHI Grant Numbers JP18H01135, JP20H01823, JP20KK0058, and~JP21H04431, JST PRESTO Grant Number JPMJPR16EA, and JST CREST Grant Number JPMJCR2014. {C.M.O.} would like to thank the financial support of CNPq (National Council for Scientific and Technological Development) grant 305383/2019-1 and FAPESP grant 2013/07375-0.
%
%
%
%
%
%

\end{document}